\documentclass[11pt]{article}
\usepackage[T1]{fontenc}
\usepackage{lmodern,amsmath,amssymb,amsthm,mathtools}
\usepackage[margin=1in]{geometry}
\usepackage{microtype,enumitem}
\usepackage{tikz}
\usepackage[hidelinks]{hyperref}
\hypersetup{pdftitle={Relative divisors and intersections of noncommutative rational localizations},pdfauthor={Mahesh Ramani}}
\setlist[enumerate]{label=\textup{(\roman*)},itemsep=2pt}
\numberwithin{equation}{section}
\newtheorem{theorem}{Theorem}[section]
\newtheorem{proposition}[theorem]{Proposition}
\newtheorem{lemma}[theorem]{Lemma}
\newtheorem{corollary}[theorem]{Corollary}
\theoremstyle{definition}
\newtheorem{definition}[theorem]{Definition}
\newtheorem{example}[theorem]{Example}
\theoremstyle{remark}
\newtheorem{remark}[theorem]{Remark}
\newcommand{\C}{\mathbb C}
\newcommand{\Q}{\mathbb Q}
\newcommand{\R}{\mathbb R}
\newcommand{\N}{\mathbb N}
\newcommand{\Z}{\mathbb Z}
\newcommand{\Mat}{\operatorname{M}}
\newcommand{\GL}{\operatorname{GL}}
\newcommand{\dom}{\operatorname{dom}}
\newcommand{\supp}{\operatorname{supp}}

\newcommand{\coker}{\operatorname{coker}}
\newcommand{\im}{\operatorname{im}}
\newcommand{\Tor}{\operatorname{Tor}}
\newcommand{\Hull}{\operatorname{Hull}}
\newcommand{\Sat}{\operatorname{Sat}}
\newcommand{\A}{\mathcal A}
\newcommand{\U}{\mathcal U}
\newcommand{\Fgroup}{\mathsf F}
\newcommand{\alg}[1]{\langle #1\rangle}
\title{Relative divisors and intersections of noncommutative rational localizations}
\author{Mahesh Ramani}
\date{}
\begin{document}
\maketitle
\begin{abstract}
The simple constituents of monic pencils determine singularity support in a free skew field. For a pencil inverse algebra that is a free ideal ring, the multiplicities of the nonbase constituents compute the scalar-regular part of Cohn's matrix divisor theory. For independent generic invertible matrices of size at least two, every entry of a nonidentity free-group word is an atom with one relative singularity component. In rank two, opposite pivots give the exact identifications among these components. Intersections of word-entry localizations are then given by intersections of their label sets. An irreducible pencil of dimension six represents a commutator entry.

For an algebra finitely generated over the free polynomial algebra at a common regular scalar point, the residual closure is a finite pencil inverse algebra determined by the generators' common domain. Every nonzero inverse-adjoining defect is infinite-dimensional. We compute how free-group quotients change intersections and give additive formulas for general quotient and tensor intersections. We prove Berenstein--Retakh's polygon conjecture using $n-2$ charts for $n\geq4$. For an unpunctured M\"obius strip with $n\geq2$ boundary marks, the upper-algebra equality follows from $n+1$ charts. For the once-marked strip, the allowed algebra is strictly smaller than its residual hull, which equals the full chart.
\end{abstract}
\tableofcontents

\section{Introduction}
Let $P=\C\alg{x_1,\ldots,x_g}$ and let $F$ be its universal skew field of fractions. A \emph{pencil inverse algebra} is a subalgebra of $F$ generated by $P$ and finitely many entries of inverses of affine pencils. We identify a universal inverse-adjoining presentation with its image in $F$ once its map into $F$ has been shown to be injective.

Cohn's divisor group is the free abelian group on stable-association classes of matrix atoms. Its coordinates count atomic factors \cite[Sections 3.2 and 7.10--7.11]{Cohn}. Free-locus and realization theorems identify singularity components with simple coefficient representations of monic pencils \cite{KV,HKV,VolDomains}. We use this identification to compute divisors over a fir $R=P\alg{L_0^{-1}}\subseteq F$. For a square matrix $A$ over $R$ invertible at a scalar point of the base domain, we obtain
\[
 \nu_R(A)=\sum_{s\notin\supp L_0}m_s(A)[s],\qquad
 \lambda_R(A)=\sum_s m_s(A).
\]
The degree is zero precisely for units and one precisely for matrix atoms. The multiplicities equal the eventual exponents of generic determinants on the base domain (Theorem~\ref{relative-theorem}).

For independent generic invertible $d\times d$ matrices $X_1,\ldots,X_r$, put
\[
 R_{r,d}=\C\alg{(X_a)_{ij},(X_a^{-1})_{ij}:a,i,j},\qquad H_w=w(X_1,\ldots,X_r).
\]
For $d\geq2$, every entry of $H_w$, $w\neq1$ in the free group $\Fgroup_r$, has relative degree one (Theorem~\ref{word-theorem}). This includes inverse letters, repetitions, commutators, and proper powers. For $d=2$, the relative classes have exactly the identifications
\begin{equation}\label{intro-pivot}
 (w,i,j)\sim(w^{-1},3-j,3-i).
\end{equation}
Writing $R_E=R_{r,2}\alg{(H_w)_{ij}^{-1}:(w,i,j)\in E}$, with $E$ a set of these labels, we prove
\begin{equation}\label{intro-intersection}
 \bigcap_\alpha R_{E_\alpha}=R_{\cap_\alpha E_\alpha}
\end{equation}
for every nonempty family. The universal presentations are faithful. A finite member supplies a finite intersection certificate (Theorem~\ref{word-intersection}).

A second part concerns an algebra $U=P\alg{f_1,\ldots,f_m}$ whose generators are regular at a common scalar point. Classical recognizable-series theory \cite{BRSeries} gives bounded observable and reachable spaces. We use these spaces to identify its left-and-right residual closure with a finite pencil inverse algebra. This closure is also the algebra of all rational functions regular on the common domain of the generators (Theorems~\ref{residual-theorem} and~\ref{envelope}). Any nonzero additive defect between $U$ and this closure is infinite-dimensional.

For a quotient group with a specified free basis, distinct labels can have the same image. An inverse can then belong to every chart image even when it is absent from the image of the original intersection.

We prove Berenstein--Retakh's polygon conjecture \cite[Conjecture 2.20]{BR}. For $n\geq4$, the fan and its flips give the equality using $n-2$ charts. We also prove the upper-algebra equality for an unpunctured M\"obius strip with $n\geq2$ boundary marks using $n+1$ charts.

\section{Pencils, domains, and singularity support}\label{sec-support}
We work over $\C$ unless a different field is specified. All algebras are unital, and all specified inverses are two-sided. A rational function is an element of $F$, equivalently an equivalence class of noncommutative rational expressions under matrix evaluation. Its domain is the union of the domains of its scalar rational representatives, at every matrix size.

Fix a scalar point $\beta\in\C^g$ and write $y_i=x_i-\beta_i$. A monic pencil at $\beta$ has the form
\[
 L=I-\sum_{i=1}^g A_i y_i,\qquad A_i\in\Mat_n(\C).
\]
Its evaluation is $I_n\otimes I_N-\sum_i A_i\otimes(X_i-\beta_i I_N)$. Write $Z(L)$ for its graded singularity locus and
\[
 R(L)=P\alg{(L^{-1})_{ab}:a,b}.
\]
A scalar-regular rational function has a minimal recognizable realization
\begin{equation}\label{realization}
 f=cL_f^{-1}b,
\end{equation}
with constant row $c$ and column $b$. Here minimal means least state dimension among representations of this form, including the constant coefficient. This convention also allows polynomial functions. We use the following established results.

\begin{theorem}[Realization and free-locus inputs]\label{inputs}
\begin{enumerate}
\item Minimal realizations exist, are unique up to constant similarity, and are reachable and observable \cite[Part I]{BRSeries}.
\item If $f$ is regular at $\beta$, then
$\dom f=\{X:L_f(X)\text{ is invertible}\}$
\cite[Theorem 3.10]{VolDomains}. Here $\dom f$ agrees with the stable extended domain. The ordinary extended domain can be larger.
\item For monic pencils, inclusion of free loci is equivalent to the coefficient-preserving homomorphism of their coefficient algebras modulo radicals in the reverse direction \cite[Theorem 3.6]{KV}.
\item If the coefficients of a nonconstant monic pencil act irreducibly, its generic determinant is irreducible at all sufficiently large matrix sizes \cite[Theorem 3.4]{HKV}.
\end{enumerate}
\end{theorem}
When an inverse matrix expression is defined at a tuple, each of its entries has a scalar rational representative defined there \cite[Lemma 3.9]{VolDomains}. Thus expressions with matrix inverses may be used in domain comparisons without choosing one scalar representative valid everywhere.

\begin{definition}
The \emph{support} $\supp L$ is the set of simultaneous-similarity classes of the nonzero irreducible coefficient tuples on the composition factors of $\C^n$ for the unital algebra generated by the $A_i$. The one-dimensional zero tuple is omitted. Multiplicities are discarded in support. For $f$, define $\supp f=\supp L_f$ using a minimal realization at the specified center.
\end{definition}
The omitted tuple gives the identity pencil. A nonzero irreducible tuple generates the full matrix algebra by Burnside's theorem. The use of a unital coefficient algebra simply records the zero tuple separately; it agrees with the nonunital formulation of Theorem~\ref{inputs}(iii).

\begin{lemma}\label{constituents}
For representatives $L_s$ of the classes in $\supp L$,
\begin{equation}\label{constituent-algebra}
 R(L)=P\alg{(L_s^{-1})_{ab}:s\in\supp L}.
\end{equation}
\end{lemma}
\begin{proof}
A composition series puts all coefficients simultaneously in block upper triangular form. The diagonal blocks of $L^{-1}$ are exactly the inverses of its diagonal pencils, so they belong to $R(L)$. Conversely triangular inversion expresses every entry of $L^{-1}$ as a finite sum of products of diagonal inverse entries and off-diagonal polynomial blocks. Constant changes of basis preserve the generated algebra.
\end{proof}

\begin{proposition}[Support membership]\label{membership}
For $f$ regular at $\beta$ and a monic pencil $L$ at $\beta$,
\begin{equation}\label{support-membership}
 f\in R(L)\quad\Longleftrightarrow\quad\supp f\subseteq\supp L.
\end{equation}
In particular, a monic pencil $K$ is invertible over $R(L)$ if and only if $\supp K\subseteq\supp L$.
\end{proposition}
\begin{proof}
If $f\in R(L)$, local scalar representatives for the inverse entries give $\dom f\supseteq\{L\text{ invertible}\}$. The domain theorem implies $Z(L_f)\subseteq Z(L)$. The free-locus inclusion theorem gives a coefficient-preserving surjection of semisimple coefficient algebras. Their simple representations are consequently included as stated. To see that all simple representations are detected by the natural composition series, use a Wedderburn splitting: the identity idempotent of an absent simple factor would act strictly upper triangularly on a faithful module and hence would be a nilpotent idempotent. Conversely, support containment and Lemma~\ref{constituents} express $L_f^{-1}$ over $R(L)$, and hence express $f$ there. Applying the same domain argument to all entries of $K^{-1}$ proves the final assertion; if they were all regular where $K$ is singular, evaluation of $KK^{-1}=I$ would be impossible.
\end{proof}

For a possibly infinite set $T$ of simple monic classes at the same center, set
\[
 R_T=P\alg{(L_s^{-1})_{ab}:s\in T}.
\]
Every element uses only finitely many $s$.

\begin{theorem}[Support intersections]\label{support-intersections}
For a nonempty indexed family $(T_\alpha)$,
\begin{equation}\label{support-intersections-eq}
 \bigcap_\alpha R_{T_\alpha}=R_{\cap_\alpha T_\alpha}.
\end{equation}
If $T_0$ is finite, that member and at most $|T_0\setminus\cap_\alpha T_\alpha|$ further members have the same intersection.
\end{theorem}
\begin{proof}
Membership of $f$ in $R_T$ is equivalent to $\supp f\subseteq T$, by Proposition~\ref{membership} and the finite-support observation. This proves the equality. For each class of $T_0$ absent from the full intersection, select one member omitting it. The selected family has the same common support.
\end{proof}

\section{Relative matrix divisors}\label{sec-divisors}
A square matrix is \emph{full} if it has no factorization through a free module of smaller rank. Two square matrices, possibly of different sizes, are \emph{stably associated} if identity stabilizations of them differ by multiplication on the left and right by invertible matrices. A \emph{matrix atom} is a full nonunit with no factorization into two square nonunits of the same size.

A left and right \emph{fir} is a ring whose left and right ideals are free of unique rank. We use the standard consequences that submodules of free modules are free, that finitely generated such submodules have finite rank, and that a fir is a Sylvester domain with a universal skew field of fractions \cite[Chapters 2, 5, and 7]{Cohn}. An embedding is \emph{honest} if it preserves fullness of matrices.

Fix a monic pencil $L_0$ and write
\[
 R=R(L_0)\subseteq F,\qquad
 \Omega_N=\{X:L_0(X)\text{ is invertible}\},\qquad B=\supp L_0.
\]
For this section assume that $R$ is a fir. Matrices under consideration are square matrices over $R$ invertible at some scalar point of $\Omega$. All pencils in an argument are recentered and normalized at that point.

\begin{lemma}[Polynomial linearization]\label{linearization}
Every square matrix $A$ over $R$ is stably associated over $R$ to a polynomial matrix. If $A$ is invertible at the chosen scalar point, it is stably associated to a monic affine pencil. These transformations preserve the algebra over $R$ generated by its inverse entries.
\end{lemma}
\begin{proof}
Induction on expressions gives
$A=D_0+CD^{-1}B'$,
where the four matrices are polynomial and $D$ is invertible over $R$ and on $\Omega$. Sums use block diagonal denominators. For a product of two such expressions, take denominator
$\left(\begin{smallmatrix}D_1&-B'_1C_2\\0&D_2\end{smallmatrix}\right)$,
row $(C_1,D_{01}C_2)$, column $(B'_1D_{02},B'_2)^T$, and polynomial term $D_{01}D_{02}$.
For $K=\left(\begin{smallmatrix}D_0&C\\-B'&D\end{smallmatrix}\right)$ one has
\[
 \begin{pmatrix}I&-CD^{-1}\\0&I\end{pmatrix}
 K\begin{pmatrix}I&0\\D^{-1}B'&I\end{pmatrix}
 =\begin{pmatrix}A&0\\0&D\end{pmatrix}.
\]
Since $D$ is a unit over $R$, this proves the first assertion. Higman's elementary linearization replaces each polynomial product by an auxiliary identity block; the identity
\[
 \begin{pmatrix}I&a_1\\0&I\end{pmatrix}
 \begin{pmatrix}a_0+a_1a_2&0\\0&I\end{pmatrix}
 \begin{pmatrix}I&0\\-a_2&I\end{pmatrix}
 =\begin{pmatrix}a_0&a_1\\-a_2&I\end{pmatrix}
\]
justifies each replacement. Constant normalization at the chosen point makes the result monic. Inverse matrices transform by the same invertible operations, which proves the last assertion.
\end{proof}

\begin{lemma}\label{honesty}
If $R(L_0)$ is a fir, its displayed embedding in $F$ is honest.
\end{lemma}
\begin{proof}
Let $A$ be full over $R$. Linearize it to a polynomial matrix $K$ by Lemma~\ref{linearization}, without using a scalar point. Stable fullness in a fir makes $K$ full over $R$. It is then full over $P$, since a smaller polynomial factorization would also be one over $R$. The universal field $F$ inverts $K$, and hence inverts $A$. This is honesty for an embedding into a skew field.
\end{proof}

\begin{lemma}[Compression]\label{compression}
Over $R$, matrix atomicity is invariant under identity stabilization and stable association.
\end{lemma}
\begin{proof}
A proper factorization stabilizes. Conversely suppose a stabilization of $A\in\Mat_n(R)$ factors as $UV$ with $U,V$ full square nonunits. Full matrices are injective, by Lemma~\ref{honesty}. For right modules the factorization gives
\[
 0\longrightarrow N\longrightarrow\coker A\longrightarrow Q\longrightarrow0,
 \qquad N\simeq\coker V,\quad Q\simeq\coker U.
\]
Both end modules are nonzero and have finite free resolutions of Euler characteristic zero. Let $I\subseteq R^n$ be the inverse image of $N$. It is generated by the columns of $A$ and finitely many lifts of generators of $N$, so it is finite free, say of rank $k$. From
$0\to AR^n\to I\to N\to0$
one gets $k-n=\chi(N)=0$. Independence of this Euler characteristic follows from Schanuel's lemma and invariant basis number. Choosing a basis of $I$ gives $A=BC$, with $B,C$ square of order $n$ and cokernels $Q,N$. They are full because their product is full, and neither is invertible. Thus $A$ has a proper factorization. The same reasoning applies after invertible changes of basis.
\end{proof}

\begin{definition}\label{relative-definition}
Choose a monic polynomial linearization $L$ of $A$. For a nonbase simple class $s\notin B$, let $m_s(A)$ be its multiplicity among the coefficient composition factors of $L$. Define
\begin{equation}\label{relative-divisor}
 \nu_R(A)=\sum_{s\notin B}m_s(A)[s],\qquad
 \lambda_R(A)=\sum_{s\notin B}m_s(A).
\end{equation}
\end{definition}

\begin{theorem}[Relative divisor theorem]\label{relative-theorem}
Definition~\ref{relative-definition} is independent of linearization. On scalar-regular square matrices over $R$ it has the following properties.
\begin{enumerate}
\item Stable association preserves $\nu_R$, and
$\nu_R(AB)=\nu_R(A)+\nu_R(B)$ whenever $A,B$ are invertible at a common scalar point. It is also additive on block sums.
\item $\lambda_R(A)=0$ if and only if $A$ is a unit; $\lambda_R(A)=1$ if and only if $A$ is a matrix atom.
\item Complete atomic factorizations exist. Every such factorization has length $\lambda_R(A)$ and the same multiset of relative simple classes.
\item Two atoms in this class are stably associated over $R$ if and only if their relative simple classes agree.
\item For all sufficiently large $N$, $\det A(\Xi_N)$ in
\[
 \mathcal O(\Omega_N)=\C[\text{entries of }\Xi_N][\det L_0(\Xi_N)^{-1}]
\]
has, up to a unit, exactly the irreducible factors indexed by $\nu_R(A)$, with exponents $m_s(A)$.
\end{enumerate}
\end{theorem}
\begin{proof}
The displayed coordinate ring is a UFD. For finitely many simple pencils, their determinants are irreducible for all large $N$ by Theorem~\ref{inputs}(iv). Distinct classes have nonassociate determinants eventually: the free-locus theorem supplies a tuple distinguishing any two free loci, and adjoining copies of the common regular scalar point gives a witness at every larger size. For a nonbase class, apply the same argument against the direct sum of the base pencils. Its determinant remains a nonunit after localization. Indeed containment of its free locus in the base locus would imply that it was a base class by Proposition~\ref{membership}.

Triangularizing $L$ factors its generic determinant as the product of its constituent determinants with their composition multiplicities. Stable association over $R$ changes determinants by a unit of $\mathcal O(\Omega_N)$, because both the change matrices and their inverses evaluate there. Unique factorization proves independence, the multiplicity assertion, and all additivity statements.

If the degree is zero, every diagonal pencil of $L$ is invertible over $R$; triangular inversion and stable association invert $A$. Units conversely have unit determinants. If the degree is at least two, cut a triangular composition series into blocks $L_1,L_2$ each containing a nonbase constituent. Then
\[
 \begin{pmatrix}L_1&C\\0&L_2\end{pmatrix}
 =\begin{pmatrix}I&0\\0&L_2\end{pmatrix}
  \begin{pmatrix}L_1&C\\0&I\end{pmatrix}
\]
is a proper factorization. Compression shows that $A$ is not an atom. If the degree is one, a proper factorization $A=BC$ would have both factors invertible at the same scalar point and of positive degree, contradicting additivity. Inductively splitting a nonatom gives a complete factorization, because positive degrees decrease. Additivity proves its asserted length and multiset.

When the degree is one, Gaussian elimination removes all base and identity blocks of $L$, leaving its unique nonbase irreducible block up to stable association. This proves the classification of atoms. The converse is divisor invariance.
\end{proof}

\begin{remark}\label{cohn-priority}
Cohn's simple $\N$-values and divisor group give factor multiplicities and their invariance \cite[Theorems 7.10.6--7.10.8]{Cohn}. Cohn also describes the surviving matrix primes and their divisors under factor-complete localization \cite[Theorem 7.11.9 and Proposition 7.11.10]{Cohn}. Theorem~\ref{relative-theorem} identifies the scalar-regular primes with relative free-locus components. Their multiplicities are the eventual exponents of the corresponding generic determinant factors. The classification up to stable association in part~(iv) concerns matrix atoms, which have degree one.
\end{remark}
Global free loci let us compare classes at different scalar centers. Recentring preserves the unital coefficient algebra and hence irreducibility. Write $C=L(\gamma)$ and $B_i=C^{-1}A_i$. Then $C^{-1}=I+\sum_i B_i(\gamma_i-\beta_i)$. In a finite-dimensional unital matrix algebra the inverse of an invertible element is polynomial in that element, so the algebras generated by the $A_i$ and the $B_i$ agree. The global free loci give the same labels at either center. Thus the divisor is intrinsic to the base algebra.

\section{Universal matrix holonomies}\label{sec-holonomy}
Let $\Gamma$ be a group and $d\geq2$. With $C=\Mat_d(\C)$ and matrix units $e_{ij}$, put
\begin{equation}\label{corner}
 B_\Gamma=C*_{\C}\C\Gamma,\qquad e=e_{11},\qquad
 R_{\Gamma,d}=eB_\Gamma e,\qquad (H_w)_{ij}=e_{1i}we_{j1}.
\end{equation}
All free products in this paper are amalgamated only over the indicated scalar field unless otherwise specified.

\begin{proposition}\label{universal-corner}
The algebra $R_{\Gamma,d}$ represents homomorphisms $\Gamma\to\GL_d(T)$ for unital $\C$-algebras $T$. In particular, for $\Gamma=\Fgroup_r$ it is the universal algebra on $r$ invertible $d\times d$ matrices.
\end{proposition}
\begin{proof}
The mutually inverse matrix-unit maps are
\[
 B_\Gamma\longrightarrow\Mat_d(eB_\Gamma e),\quad
 b\longmapsto(e_{1i}be_{j1})_{ij},\qquad
 (b_{ij})\longmapsto\sum_{ij}e_{i1}b_{ij}e_{1j}.
\]
Inserting $1=\sum_k e_{k1}e_{1k}$ gives $H_wH_v=H_{wv}$. A representation of $\Gamma$ in $\GL_d(T)$ and the fixed matrix units define a map $B_\Gamma\to\Mat_d(T)$, whose corner is a map to $T$. Conversely the matrix-unit isomorphism recovers that representation from a corner map. These constructions are inverse.
\end{proof}

Let $P_{r,d}$ be free on the entries of $r$ generic matrices $X_a$, and $F_{r,d}$ its free skew field.

\begin{lemma}[Faithfulness]\label{matrix-magnus}
The universal algebra $R_{\Fgroup_r,d}$ embeds into $F_{r,d}$ by $H_{z_a}\mapsto X_a$. Its image is $R_{r,d}$.
\end{lemma}
\begin{proof}
The universal presentation is $P_{r,d}\alg{X_1^{-1},\ldots,X_r^{-1}}$. To prove injectivity, use the Magnus--Fox embedding
\[
 \C\Fgroup_r\hookrightarrow\C\alg{\!\alg{u_1,\ldots,u_r}\!},\qquad z_a\mapsto1+u_a;
\]
see \cite[Theorem 2.11]{AD}. Split the group algebra into scalars and its augmentation ideal $I$, and split $C=\C1\oplus C^\circ$. Free-product normal forms are alternating tensor words from $I$ and $C^\circ$. The expansion of every $I$-slot has zero constant coefficient. Consequently successive $C^\circ$ slots remain separated by nonempty $u$-words in the $u$-degree completion of $C*\C\alg{u_1,\ldots,u_r}$. Different slot patterns remain distinct. Within one pattern, finite independent families of series are separated by finitely many coefficient functionals; tensoring these functionals proves injectivity on every finite tensor sum. Thus $C*\C\Fgroup_r$ embeds in this completion.

The corner of the uncompleted algebra is free on the $d^2r$ entries $e_{1i}u_ae_{j1}$, by the arbitrary-matrix universal property. Its completion is the ordinary free power-series ring in those entries. The corner expansion therefore embeds $R_{\Fgroup_r,d}$ and sends $H_{z_a}$ to $I+U_a$. It is the rational expansion at $X_a=I$. A zero rational function regular there has zero series: nilpotent chain evaluations extract each coefficient and make inverse expansions finite. Hence the map to $F_{r,d}$ is injective. Its image has the stated generators.
\end{proof}

We recall the particular localization results used below. A set $\Sigma$ of full square matrices is \emph{factor-inverting} if square factors of a member become invertible in the universal localization. One may pass to its multiplicative completion, consisting of invertible changes of block upper triangular matrices with diagonal blocks in $\Sigma\cup\{I\}$. Over a semifir this preserves the factor-inverting property \cite[Lemma 7.11.1]{Cohn}. Cohn's applicable statements are
\begin{equation}\label{cohn-localization}
 \begin{gathered}
 R\text{ a fir},\ \Sigma\text{ factor-inverting}\ \Longrightarrow\ R_\Sigma\text{ a fir},\\
 f:R\to D\text{ honest and }\Sigma\text{-inverting}
 \ \Longrightarrow\ R_\Sigma\to D\text{ injective}.
 \end{gathered}
\end{equation}
These are \cite[Theorem 7.11.7 and Proposition 7.5.9(ii)]{Cohn}. In each application we use the multiplicative completion of the initial set.

\begin{proposition}\label{holonomy-fir}
$R_{r,d}$ is a fir, and its embedding in $F_{r,d}$ is honest.
\end{proposition}
\begin{proof}
Each $X_a$, normalized at the identity, is a matrix atom over $P_{r,d}$. Indeed its evaluation is an ordinary generic matrix of size $dN$, whose determinant is irreducible at every $N$. The monic-pencil atom criterion \cite[Proposition 4.5(2)]{HKV} applies at a sufficiently large $N$. If $AB=X_a$ with square polynomial factors, one factor is a polynomial unit and both factors become units when $X_a$ is inverted. The set of these atoms is therefore factor-inverting, as is its multiplicative completion. Apply \eqref{cohn-localization} and identify the universal localization with the faithful image from Lemma~\ref{matrix-magnus}. Lemma~\ref{honesty} gives honesty.
\end{proof}

\begin{lemma}[Corner degree]\label{corner-degree}
For every group $\Gamma$, the corner in \eqref{corner} is a domain with only scalar units. Every entry $(H_w)_{ij}$, $w\neq1$, is an atom.
\end{lemma}
\begin{proof}
Filter free-product normal forms by the number of nonidentity group blocks. The degree-zero corner is $eCe=\C e$. For positive degrees $n,m$, write leading corner tensors as
$\sum_s U_s\otimes c_s$ and $\sum_t r_t\otimes V_t$, with independent outer families $U_s,V_t$ and nonzero $c_s\in Ce$, $r_t\in eC$. The leading degree-$n+m$ product is
\[
 \sum_{s,t} U_s\otimes[c_sr_t]\otimes V_t,
\]
where the middle class lies in $C/\C1$. The product of a nonzero column and nonzero row has rank one and is not scalar for $d\geq2$. Each middle class is therefore nonzero, and independence of the outer families prevents cancellation. Degree is additive on nonzero products. Every $(H_w)_{ij}$, $w\neq1$, has degree one: its nonzero group block is $w$ with nonzero outer matrix-unit factors. In a factorization one factor must consequently have degree zero and be scalar. Additivity also proves the domain and unit assertions.
\end{proof}
\section{Arbitrary words and their relative components}\label{sec-words}
\begin{theorem}\label{word-theorem}
Let $d\geq2$, $r\geq1$, and $w\neq1$ in $\Fgroup_r$. Every $h=(H_w)_{ij}$ has
\begin{equation}\label{word-degree}
 \lambda_{R_{r,d}}(h)=1.
\end{equation}
The algebra $R_{r,d}[h^{-1}]$ has precisely one additional simple-pencil class. For all sufficiently large $N$, the hypersurface
\[
 \{X\in\GL_{dN}(\C)^r:\det h(X)=0\}
\]
is nonempty, reduced, and irreducible. If $w$ has reduced length $\ell$, the additional class has an irreducible pencil of dimension at most $d\ell+1$.
\end{theorem}
\begin{proof}
Scalar atomicity follows from Lemma~\ref{corner-degree} and faithfulness; the fir property is Proposition~\ref{holonomy-fir}. Every entry admits a scalar point where it is invertible. For diagonal entries take all $X_a=I$. For off-diagonal entries use the faithful projective free-group representation constructed in Appendix~\ref{appendix-pingpong}. Every nonidentity word is nonscalar, and a constant conjugation makes any specified off-diagonal entry nonzero. Embedding the $2\times2$ example into size $d$ and conjugating proves the general case.

The relative divisor theorem now gives \eqref{word-degree} and the hypersurface conclusion. To bound the pencil dimension, write the word as $\ell$ signed letters and assign
\[
 (D_k,C_k)=(I,X_a)\text{ for }z_a,\qquad
 (D_k,C_k)=(X_a,I)\text{ for }z_a^{-1}.
\]
Let $D$ be block upper bidiagonal of size $d\ell$ with diagonal $D_k$ and superdiagonal $-C_k$. Let $b$ have final block $C_\ell e_j$ and other blocks zero, and let $c$ select the $i$th coordinate of its first block. Back substitution gives $h=cD^{-1}b$. The affine Schur pencil
\begin{equation}\label{word-pencil}
 A_w=\begin{pmatrix}0&c\\-b&D\end{pmatrix}
\end{equation}
is stably associated to $h$ over the base. It is invertible at the selected scalar point, so its unique nonbase constituent has dimension at most $d\ell+1$.
\end{proof}
For a finite list a common scalar point is obtained by intersecting nonempty Zariski open sets. For all word entries simultaneously, one can choose the entry coordinates algebraically independent over $\Q$: each required nonvanishing condition is a nonzero rational function with rational coefficients. This fixes a common center for the support calculations below.

\section{Rank-two separation and intersections}\label{sec-intersections}
\begin{lemma}[Opposite pivots]\label{opposite}
For an invertible $2\times2$ matrix $H$ over a unital ring, adjoining an inverse to $H_{ij}$ is universally equivalent to adjoining an inverse to $(H^{-1})_{3-j,3-i}$.
\end{lemma}
\begin{proof}
Permute rows and columns to reduce to $a=H_{11}$. If $a$ is a unit, Gaussian factorization shows that the Schur complement $s=d-ca^{-1}b$ is a unit and $(H^{-1})_{22}=s^{-1}$, where $H=\left(\begin{smallmatrix}a&b\\c&d\end{smallmatrix}\right)$. Apply this implication to $H^{-1}$ for the converse. The resulting universal maps fix the base and are inverse.
\end{proof}

\begin{lemma}[Plane extension]\label{plane-extension}
Let $\dim V=2N$, $\dim E=t$, $T\in\GL(V)$, and $E\cap TE=0$. An $N$-plane $W\supseteq E$ with $W\cap TW=0$ exists if either
\begin{enumerate}
\item $N\geq t+2$ and every eigenspace of $T$ has dimension at most $N-t$; or
\item $T$ is an involution whose two eigenspaces have dimension $N$.
\end{enumerate}
The set of these planes is nonempty and open in the Grassmannian of planes containing $E$.
\end{lemma}
\begin{proof}
That Grassmannian is irreducible of dimension $N(N-t)$. A bad plane contains $x,Tx$ for a line $[x]$. Eigenlines have parameter dimension at most $N-t-1$, and requiring their additional line has codimension $N$. If $x,Tx$ are independent and their two-plane misses $E$, containment has codimension $2N$, while $[x]$ varies in dimension at most $2N-1$. If the two-plane meets $E$ in a line, then $(\alpha I+\beta T)x\in E$ for some $[\alpha:\beta]\in\mathbb P^1$. At nonsingular parameters this family of lines has dimension at most $t$; at each of the finitely many singular parameters it has dimension at most
$t+\dim\ker(\alpha I+\beta T)-1\leq N-1$.
Containment has codimension $N$. A two-plane contained in $E$ is excluded by $E\cap TE=0$. Each incidence image has smaller dimension, proving (i). For (ii), write $V=V_+\oplus V_-$. The projections of $E$ are injective; thus $E$ is a partial graph isomorphism. Extend it to an isomorphism $V_+\to V_-$. Its graph and its negative graph are complementary. Openness follows from a determinant condition.
\end{proof}

\begin{theorem}[Matrix-evaluation separation]\label{separation}
Let $\Gamma$ be residually finite. Fix $w\neq1$ and $(i,j)\in\{1,2\}^2$. For any finite list of nonidentity-word entries excluding
$(w,i,j)$ and $(w^{-1},3-j,3-i)$, there is a homomorphism
$R_{\Gamma,2}\to\Mat_N(\C)$
making the selected entry singular and all listed entries invertible.
\end{theorem}
\begin{proof}
Simultaneous coordinate interchange reduces to selected type $(1,1)$ or $(1,2)$. For $(1,1)$ set $\epsilon_1=0$, $\epsilon_2=1$. For each listed triple $(v,k,l)$ put
$a=w^{-\epsilon_{3-k}}vw^{\epsilon_l}$.
If $a=1$, the assumption $v\neq1$ gives exactly an excluded triple. Choose a finite quotient $G$ preserving these nonidentities and $w\neq1$. On two copies of the regular representation, of dimension $2N=2|G|$, every nonidentity element has eigenspaces of dimension at most $N$. Pair eigenvectors of different eigenvalues to obtain an $N$-plane complementary to its translate. A pairing exists by repeatedly removing vectors from the two largest eigenspaces; their sizes remain at most half the remaining total. Consequently each required complementarity condition is nonempty open. Choose $U$ satisfying all of them and put $W=\rho(w)U$. In $V=U\oplus W$ the selected block is zero. A listed block is invertible precisely when
\[
 \rho(v)\rho(w)^{\epsilon_l}U\oplus\rho(w)^{\epsilon_{3-k}}U=V,
\]
which is the imposed condition for $a$.

For selected type $(1,2)$, choose a finite quotient retaining the listed words and $w$ as nonidentity and retaining the four distinct elements $1,w,v,vw$ for every listed $(v,1,2)$. These four are distinct because $v\notin\{1,w,w^{-1}\}$. Take twelve regular copies, so $2N=12|G|$, and in one copy put $u=\delta_1$, $E=\operatorname{span}(u,\rho(w)u)$. For every listed $(v,1,2)$, the four basis vectors give $E\cap\rho(v)E=0$. If $v$ has order two, use Lemma~\ref{plane-extension}(ii); otherwise its eigenspaces have dimension at most $2N/3\leq N-2$, so use (i) with $t=2$. Intersect these nonempty opens to choose $W\supseteq E$ with $W\cap\rho(v)W=0$ for each listed type $(1,2)$. The selected block is singular because $\rho(w)u\in W\cap\rho(w)W$.

Choose an $N$-plane $U$ satisfying $U\cap W=0$ and, for the other listed types,
\[
 \begin{array}{c|c}
 (1,1)&U\cap\rho(v)^{-1}W=0\\
 (2,2)&U\cap\rho(v)W=0\\
 (2,1)&U\cap\rho(v)U=0.
 \end{array}
\]
The first two are open complementarity conditions with a fixed $N$-plane; each is nonempty. The third is nonempty by the eigenvector-pairing argument. Their finite intersection in the full irreducible Grassmannian is nonempty. In $U\oplus W$ these are exactly the remaining block invertibility conditions. The universal property in Proposition~\ref{universal-corner} gives the required corner map.
\end{proof}
Free groups are residually finite. For completeness, a finite neighborhood in their Cayley tree containing a reduced-word path carries partial permutations for the generators. Complete each partial injection to a permutation of that finite vertex set. The path still has distinct endpoints. Products of these finite permutation quotients separate any finite list of nonidentity words.

\begin{theorem}[Classification and arbitrary intersections]\label{word-intersection}
For $R=R_{r,2}$, the relative simple classes of all nonidentity-word entries are indexed exactly by triples modulo \eqref{intro-pivot}. For sets $E$ of these labels, the algebras $R_E$ have faithful universal inverse-adjoining presentations, and for any nonempty family
\begin{equation}\label{word-intersection-eq}
 \bigcap_\alpha R_{E_\alpha}=R_{\cap_\alpha E_\alpha}.
\end{equation}
If $E_0$ is finite, it and at most $|E_0\setminus\cap E_\alpha|$ further members suffice. This bound is sharp for unrestricted families.
\end{theorem}
\begin{proof}
Every entry supplies one new class by Theorem~\ref{word-theorem}. Opposite pivots give the same inverse algebra by Lemma~\ref{opposite}, hence the same class. If two other triples gave the same class, support membership would make their inverse algebras equal. Theorem~\ref{separation} gives an evaluation where one entry is singular and the other and the base holonomies are invertible. The latter inverse algebra is defined there and cannot contain the inverse of the singular entry. This contradicts equality. Thus the list of classes is exact.

Each entry is an atom of the honest fir $R$. The set to be inverted is factor-inverting: in a square factorization of a scalar atom one scalar factor is a unit. Its multiplicative completion remains factor-inverting. Proposition 7.5.9(ii) of \cite{Cohn}, applied to $R\hookrightarrow F_{r,2}$, proves injectivity of the universal localization. For infinite sets, the same argument, or the filtered union of finite localizations, applies. The relative support of $R_E$ is therefore exactly $E$ in addition to the common base support. Theorem~\ref{support-intersections} proves the intersection and the bound.

For sharpness choose distinct labels $e_1,\ldots,e_q,f$ and set
$E_0=\{e_1,\ldots,e_q\}$,
$E_i=(E_0\setminus\{e_i\})\cup\{f\}$.
Every one of these $q+1$ members is necessary: omitting $E_0$ retains $f$, and omitting $E_i$ retains $e_i$.
\end{proof}
Theorem~\ref{word-theorem} applies in every rank $d\geq2$ to all nonidentity free-group words, including words with trivial abelianization. The classification and intersection formula in Theorem~\ref{word-intersection} are for rank two.

\section{An explicit commutator pencil}\label{sec-commutator}
Write
\[
 X=\begin{pmatrix}a&b\\c&d\end{pmatrix},\qquad
 Y=\begin{pmatrix}e&f\\g&h\end{pmatrix},\qquad
 \eta=(XYX^{-1}Y^{-1})_{11}.
\]
Consider the affine pencil, monic at $(X,Y)=(I,I)$,
\begin{equation}\label{commutator-pencil}
 L=\begin{pmatrix}
 a&b&0&0&0&0\\
 0&1&c-g&d-h&g&h-1\\
 a-1&b&e&f&0&0\\
 0&0&c&d&g&h-1\\
 a-1&b&e-a&f-b&1&0\\
 0&0&0&0&g&h
 \end{pmatrix}.
\end{equation}

\begin{theorem}\label{commutator-theorem}
The coefficient matrices of $L$ generate $\Mat_6(\C)$. Over $R_{2,2}$,
\[
 R_{2,2}[\eta^{-1}]=R_{2,2}\alg{L^{-1}}.
\]
At every matrix size, on the domain where $X,Y$ are invertible,
\begin{equation}\label{commutator-det}
 \det L(X,Y)=\det X\det Y\det\eta(X,Y).
\end{equation}
Thus $L$ represents the unique additional relative class of $\eta$, with multiplicity one.
\end{theorem}
\begin{proof}
Set
\[
 D=\begin{pmatrix}I&-X&0&0\\0&I&-Y&0\\0&0&X&-I\\0&0&0&Y\end{pmatrix},\quad
 b_0=\begin{pmatrix}0\\0\\0\\e_1\end{pmatrix},\quad
 c_0=(e_1^T,0,0,0),\quad A=\begin{pmatrix}0&c_0\\-b_0&D\end{pmatrix}.
\]
The four blocks of $D$ have size two. Back substitution gives $c_0D^{-1}b_0=\eta$. With $A_0=A(I,I)$, coefficient-by-coefficient multiplication gives
\begin{equation}\label{commutator-reduction}
 A_0^{-1}A=\begin{pmatrix}I_3&B_0\\0&L\end{pmatrix},\qquad
 B_0=\begin{pmatrix}
 a-1&b&e-a&f-b&1-e&-f\\
 0&0&0&0&0&0\\
 -c&1-d&c-g&d-h&g&h-1
 \end{pmatrix}.
\end{equation}
This is an affine polynomial identity. Its verification requires only the constant coefficient and the eight displayed variable coefficients. Schur elimination of $D$, which is invertible over the base, identifies the inverse-entry algebra of $A$ over the base with that of $\eta$. Equation~\eqref{commutator-reduction} identifies it with that of $L$. Since $\det A_0=1$ and $\det D=\det X\det Y$, the ordinary Schur determinant identity after evaluation proves \eqref{commutator-det}.

Write
$L=I+Q_a(a-1)+Q_b b+Q_c c+Q_d(d-1)+Q_e(e-1)+Q_f f+Q_g g+Q_h(h-1)$.
Let $v=e_3+e_5$, $u=e_3^T$, so that $Q_e=vu$. The following column and row calculations give a complete rank-one certificate:
\[
\begin{array}{c|c@{\qquad}c|c}
 A&Av&B&uB\\[2pt]\hline
 I&e_3+e_5&I&e_3^T\\[2pt]
 Q_a&-e_5&Q_a&e_1^T\\[2pt]
 Q_c&e_2+e_4&Q_b&e_2^T\\[2pt]
 Q_g&e_4+e_6&Q_f&e_4^T\\[2pt]
 Q_gQ_a&-e_2-e_4-e_6&Q_bQ_g&-e_3^T+e_5^T\\[2pt]
 Q_bQ_c&e_1+e_3&Q_bQ_h&-e_4^T+e_6^T
\end{array}
\]
Both lists are bases. Hence the thirty-six products $AQ_eB=(Av)(uB)$ span $\Mat_6(\C)$. The pencil is irreducible. Its class differs from the base classes, both of dimension two; its eventual irreducible determinant remains nonunit on the base. Equation~\eqref{commutator-det} gives multiplicity one and proves the last assertion.
\end{proof}
The supplementary file \texttt{verify\_commutator\_pencil.py} verifies \eqref{commutator-reduction}, both six-vector bases, and the rank of the thirty-six products using exact rational arithmetic. The proof above is independent of executing it.

\begin{example}
With $R=R_{2,2}$, Theorem~\ref{word-intersection} gives
\[
 R[\eta^{-1},a^{-1}]\cap R[\eta^{-1},e^{-1}]=R[\eta^{-1}].
\]
The three labels are distinct. In particular the commutator, whose abelianization is zero, is covered by the arbitrary-word theorem.
\end{example}

\section{Residual completion and saturation}\label{sec-residual}
Fix a scalar center $\beta$ and $y_i=x_i-\beta_i$. If
$f=\sum_w f_w y_w$ is a rational series regular there, define the coefficient residuals by
\begin{equation}\label{residual-def}
 (\ell_u f)_w=f_{uw},\qquad (r_v f)_w=f_{wv}.
\end{equation}
For letters these are $\ell_i,r_i$. Empty words give the identity operator. For $u=i_1\cdots i_k$ our convention is $\ell_u=\ell_{i_k}\cdots\ell_{i_1}$ and $r_u=r_{i_1}\cdots r_{i_k}$; operator products act from right to left. Left and right residuals commute. Residuals preserve rationality, as follows also from \eqref{residual-realization} below.

\begin{lemma}\label{residual-products}
For scalar rational series,
\begin{align}\label{residual-product}
 \ell_i(fg)&=(\ell_i f)g+f(\beta)\ell_i g,&
 r_i(fg)&=f(r_i g)+(r_i f)g(\beta).
\end{align}
If $K=(I-\sum_i A_i y_i)^{-1}$, then
\[
 \ell_i K=A_iK,\qquad r_iK=KA_i.
\]
In particular $P\alg{K_{ab}}$ is residual-stable.
\end{lemma}
\begin{proof}
In a product coefficient whose word begins with $i$, either the first factor contributes a nonempty prefix or it contributes its constant coefficient. This proves the left rule; the suffix argument proves the right rule. Expansion $K=\sum_w A_wy_w$ proves the matrix formulas. The product rules applied to the algebra generators prove stability.
\end{proof}

\begin{definition}
For $U=P\alg{f_1,\ldots,f_m}$ regular at $\beta$, let $\Hull_\beta(U)$ be the smallest subalgebra containing $U$ and stable under every left and right coefficient residual at $\beta$.
\end{definition}

\begin{theorem}[Finite residual completion]\label{residual-theorem}
Let $f_j=c_jL_j^{-1}b_j$ be minimal realizations of dimensions $d_j$. Then
\begin{equation}\label{residual-hull}
 \Hull_\beta(U)=P\alg{(L_j^{-1})_{ab}:j,a,b}.
\end{equation}
It is generated by $P$ and the finite set
\begin{equation}\label{bounded-residuals}
 \ell_u r_v f_j,\qquad |u|,|v|\leq d_j-1.
\end{equation}
In particular $U$ is residual-stable if and only if the first left and right residuals of its stated generators lie in $U$. A finitely generated algebra containing $P$ and regular at a common scalar point is a pencil inverse algebra if and only if it is residual-stable.
\end{theorem}
\begin{proof}
For one generator $f=cKb$ of minimal dimension $d$, reachability and observability give
\[
 \operatorname{span}\{A_vb\}=\C^d,\qquad
 \operatorname{span}\{cA_u\}=(\C^d)^*.
\]
Lengths at most $d-1$ suffice. The successive spans begin with a nonzero vector or row, and once a stage stops increasing it is invariant and never increases again. The coefficient convention gives
\begin{equation}\label{residual-realization}
 \ell_u r_v f=cA_uKA_vb.
\end{equation}
Choose bases among those rows and columns with bounded lengths. Every entry of $K$ is a constant linear combination of the corresponding double residuals; conversely every such residual is a constant linear combination of inverse entries. Lemma~\ref{residual-products} shows that the algebra generated by these inverse entries and $P$ is residual-stable. This proves \eqref{residual-hull}--\eqref{bounded-residuals}, and then handles all generators by taking their union. Constant generators cause no difficulty and may be omitted.

If the first residuals of the stated generators lie in $U$, the product rules make the whole algebra residual-stable. Conversely closure includes those residuals. Finally a pencil inverse algebra is residual-stable by the lemma, and every algebra equal to its hull has the finite pencil presentation \eqref{residual-hull}.
\end{proof}
The first-residual condition is a finite collection of subalgebra membership questions.

\begin{theorem}[Domain envelope]\label{envelope}
For $U=P\alg{f_1,\ldots,f_m}$ as above put $D_U=\bigcap_j\dom f_j$. Then
\begin{equation}\label{domain-envelope}
 \Hull_\beta(U)=\{f\in F:D_U\subseteq\dom f\}.
\end{equation}
Consequently the hull is independent of the common scalar expansion center and of the finite algebra generating set. Polynomial automorphisms of $P$ transport the hull covariantly.
\end{theorem}
\begin{proof}
Let $L=\bigoplus_j L_j$. The domain theorem gives $D_U=\{L\text{ invertible}\}$. Every member of $R(L)$ is regular there. Conversely a function regular on $D_U$ is regular at $\beta$, and its minimal pencil satisfies $Z(L_f)\subseteq Z(L)$. Proposition~\ref{membership} places it in $R(L)$; use Theorem~\ref{residual-theorem}. Polynomial combinations of generators are regular on their common domain. Expressing two finite generating sets in terms of one another therefore proves that their common domains agree. A polynomial automorphism and its inverse give mutually inverse maps on matrix tuples at every size; composition transports \eqref{domain-envelope}.
\end{proof}

For $R=R(L_0)$ and $U=R\alg{f_1,\ldots,f_m}$, assume
$D_{U,R}=\Omega\cap\bigcap_j\dom f_j$
contains a scalar point and define
\[
 \Sat_R(U)=\{f\in F:D_{U,R}\subseteq\dom f\}.
\]
The inverse entries of $L_0$ collectively have domain exactly $\Omega$, by evaluating its inverse identities. Thus $\Sat_R(U)$ is the same finite hull obtained by including those entries among the generators over $P$.

\begin{proposition}[Coordinate covariance]\label{covariance}
Suppose a rational automorphism $\alpha$ of $F$ is represented on every $\Omega_N$ by a bijection whose coordinate functions and inverse are everywhere regular there, and $\alpha R=R$. Then
\[
 \alpha(\Sat_R(U))=\Sat_R(\alpha U).
\]
This applies to every free-group-basis transformation of the independent holonomies over $R_{r,d}$.
\end{proposition}
\begin{proof}
Composition with the regular coordinate map transports regularity on $\Omega$; composition with its regular inverse proves the reverse implication. It therefore transports the domain envelope. A free-group automorphism substitutes matrix words for the $X_a$. These substitutions and their inverse are regular on $\GL_{dN}^r$. The universal property gives inverse automorphisms of $R_{r,d}$, and honesty extends them to the free skew field \cite[Corollary 7.5.16]{Cohn}. All hypotheses are satisfied.
\end{proof}

\begin{example}[Failure of unrestricted birational invariance]\label{birational-example}
In $F$ on $a,b$ put $P=\C\alg{a,b}$ and $U=P\alg{ab^{-1}}$, centered at $(0,1)$. Since $\ell_a(ab^{-1})=b^{-1}$,
\[
 \Hull_P(U)=P[b^{-1}].
\]
In the rational coordinates $t=ab^{-1},b$, however, $U=\C\alg{t,b}=P'$, so $\Hull_{P'}(U)=U$. These coordinates are free: $ab^{-1},b$ form a free-group basis in the containing free group algebra, and the inverse substitution is $a=tb$. The classes $b^{-n}$ are independent modulo $U$, as seen by the map $a\mapsto0$, $b\mapsto z$ into $\C[z^{\pm1}]$, which sends $U$ into $\C[z]$.
\end{example}

\begin{theorem}[Infinite additive defects]\label{infinite-defect}
Let $U\subseteq H$ be unital algebras over a field, and suppose $H$ is generated over $U$ by inverse entries of square matrices over $U$. Then $H=U$ or $\dim(H/U)=\infty$. More precisely, if $L\in\Mat_n(U)$ has inverse over $H$ but not over $U$, and $F_m\subseteq H/U$ is spanned by entries of $L^{-1},\ldots,L^{-m}$, then
\begin{equation}\label{defect-bound}
 \dim F_m\geq\lceil m/n^2\rceil.
\end{equation}
In particular a nonzero residual-saturation defect is infinite-dimensional.
\end{theorem}
\begin{proof}
The matrix classes $[L^{-1}],\ldots,[L^{-m}]$ are linearly independent modulo $\Mat_n(U)$. Otherwise take a relation with last nonzero coefficient $c_k$ and multiply on the left by $L^{k-1}$. All earlier powers become nonnegative powers over $U$, so $c_kL^{-1}\in\Mat_n(U)$, a contradiction. These $m$ independent matrix classes lie in $\Mat_n(F_m)$, of dimension $n^2\dim F_m$, proving \eqref{defect-bound}. If every adjoined inverse matrix were already over $U$, the generated algebra would be $U$; otherwise this argument applies. The inverse realization pencils of a residual completion are matrices over $P\subseteq U$.
\end{proof}
We view $H/U$ as a $U$-bimodule and take its dimension over the ground field.

\section{Constrained holonomies and quotient descent}\label{sec-descent}
\subsection{Explicitly free quotient groups}
Let $\pi:\Fgroup_r\twoheadrightarrow\Fgroup_s$, $s\geq1$, exhibit a quotient group with a specified free basis. Let $J\subseteq R_{r,2}$ be generated by the entries of $H_k-I$ for defining relators $k$. Proposition~\ref{universal-corner} gives
\begin{equation}\label{free-quotient}
 R_{r,2}/J\simeq R_{s,2}.
\end{equation}
For every unital target algebra, representations of the quotient correspond to representations that kill the relators. The specified quotient coordinates give inverse presentation maps. Lemma~\ref{matrix-magnus} gives the faithful free-field model for the new free basis.

Call a word-entry label admissible if its word has nonidentity image under $\pi$, and write $\pi_*$ for the induced map on labels modulo opposite pivots. Identity-word diagonal entries are units and may be discarded. Identity-word off-diagonal entries are zero, so a faithful inverse-adjoining quotient cannot invert them.

\begin{theorem}[Support collisions]\label{quotient-collisions}
For a nonempty family of sets of admissible labels $(E_\alpha)$ put
\[
 S_0=\pi_*\bigl(\bigcap_\alpha E_\alpha\bigr),\qquad
 S_1=\bigcap_\alpha\pi_*(E_\alpha).
\]
In the common quotient localization,
\begin{equation}\label{quotient-intersections}
 q\bigl(\bigcap_\alpha(R_{r,2})_{E_\alpha}\bigr)=(R_{s,2})_{S_0},\qquad
 \bigcap_\alpha q\bigl((R_{r,2})_{E_\alpha}\bigr)=(R_{s,2})_{S_1}.
\end{equation}
The quotient preserves this intersection if and only if $S_0=S_1$. The additional common algebra is generated over $(R_{s,2})_{S_0}$ by the inverse entries with labels in $S_1\setminus S_0$. If this difference is nonempty, its additive defect is infinite-dimensional.
\end{theorem}
\begin{proof}
Let $E=\bigcup_\alpha E_\alpha$. Universality of inverse adjunction and \eqref{free-quotient} give
\[
 (R_{r,2})_E/J(R_{r,2})_E\simeq(R_{s,2})_{\pi_*E}.
\]
Here $J(R_{r,2})_E$ means the two-sided ideal generated by $J$. Both localization presentations are faithful by Theorem~\ref{word-intersection}; admissibility prevents zero denominators in the target. Every individual chart image is exactly $(R_{s,2})_{\pi_*E_\alpha}$. Apply Theorem~\ref{word-intersection} before and after quotienting to obtain \eqref{quotient-intersections}. Distinct labels are distinct support classes, proving the equality criterion and the generation assertion. For a missing label $e$, its entry $a$ lies in the smaller algebra, while $a^{-1}$ does not. The scalar case of Theorem~\ref{infinite-defect}, or its proof applied to the powers $a^{-n}$, proves the last assertion.
\end{proof}

\begin{example}[Dependent holonomies]\label{dependent-example}
Impose $H_2=H_1$. The quotient is $R_{1,2}$. Before quotienting, for $a=(H_1)_{11}$ and $b=(H_2)_{11}$,
\[
 R_{2,2}[a^{-1}]\cap R_{2,2}[b^{-1}]=R_{2,2}.
\]
Both images after quotienting are $R_{1,2}[h_{11}^{-1}]$. The defect is therefore
$R_{1,2}[h_{11}^{-1}]/R_{1,2}$,
with independent classes $h_{11}^{-n}$. A direct witness is
\[
 H=\begin{pmatrix}t&1\\1&0\end{pmatrix},\qquad
 H^{-1}=\begin{pmatrix}0&1\\1&-t\end{pmatrix}.
\]
The base evaluates in $\Mat_2(\C[t])$, while the extra inverse evaluates to $t^{-1}$. More generally, a primitive relator eliminates one free-group-basis generator. The relation $H_1\cdots H_r=I$ eliminates the last holonomy and yields $R_{r-1,2}$. For three factors, the theorem thus applies with free coordinates $H_1,H_2$ and $H_3=(H_1H_2)^{-1}$, including every commutator word in the remaining coordinates. This relation alone includes no surface decoration constraints.
\end{example}

\subsection{An involution pivot chart}
\begin{theorem}\label{involution}
The universal algebra of $H=\left(\begin{smallmatrix}a&b\\c&d\end{smallmatrix}\right)$ subject to $H^2=I$ and $b$ invertible is
\begin{equation}\label{involution-chart}
 B\simeq\C\alg{a,b,b^{-1}},\qquad
 H=\begin{pmatrix}a&b\\b^{-1}(1-a^2)&-b^{-1}ab\end{pmatrix}.
\end{equation}
It embeds faithfully in the free skew field on $a,b$. Relative to $B$, the entry $c$ has exactly two distinct components, represented by $1-a$ and $1+a$, each with multiplicity one. At every matrix size these are two reduced irreducible hypersurfaces on $\Mat_N(\C)\times\GL_N(\C)$. Moreover
\begin{equation}\label{involution-intersection}
 B[c^{-1}]=B[(1-a)^{-1},(1+a)^{-1}],\qquad
 B[a^{-1}]\cap B[c^{-1}]=B.
\end{equation}
\end{theorem}
\begin{proof}
The equations $ab+bd=0$ and $a^2+bc=1$ force the formulas for $d,c$. Conversely
\[
 ca+dc=b^{-1}((1-a^2)a-a(1-a^2))=0,\qquad cb+d^2=1,
\]
so the remaining equations hold. These substitutions are inverse universal presentation maps. The scalar localization $\C\alg{a,b}[b^{-1}]$ embeds in the free group algebra on $a,b$ and in its free skew field, by reduced Laurent-word normal forms and the Magnus embedding. The same formulas identify every matrix evaluation chart with $\Mat_N\times\GL_N$.

Since $b$ is a unit, inverting $c$ is equivalent to inverting $1-a^2$. Its two commuting factors $1-a$ and $1+a$ then become units: a right and a left inverse are obtained by multiplying $(1-a^2)^{-1}$ by the other factor, and agree. The one-variable pencils $1-a$ and $1+a$, normalized at $a=0$, have distinct scalar coefficient classes, different also from the class of $b$ at $b=1$. Their generic determinants are translates of the generic determinant, irreducible and nonassociate at every size. They remain so on inverting $\det b$, and
\[
 \det c=\det(b)^{-1}\det(I-A)\det(I+A).
\]
This proves the component and multiplicity assertions. The class of $a$ (center at any scalar point avoiding $0,1,-1$) is distinct from both; support intersection proves the second equality in \eqref{involution-intersection}.
\end{proof}
The corner-to-chart map sends the atom $c$ of Lemma~\ref{corner-degree} to $b^{-1}(1-a)(1+a)$. Both $1-a$ and $1+a$ are nonunits in the chart.

\subsection{General additive descent formulas}
\begin{theorem}[Quotient intersections]\label{general-quotient}
Let $q:D\to D/J$ be a quotient of algebras and let $A_1,\ldots,A_m$ be subalgebras, $m\geq2$. Set $C=\bigcap_i A_i$, $K=\bigcap_i q(A_i)$, and define
\[
 \delta:\bigoplus_i A_i\to D^{m-1},\qquad
 \delta(a_1,\ldots,a_m)=(a_2-a_1,\ldots,a_m-a_1).
\]
Writing $\delta_J$ for its restriction to $\bigoplus_i(A_i\cap J)$, there is an exact sequence of vector spaces
\begin{equation}\label{descent-exact}
 0\longrightarrow q(C)\longrightarrow K\xrightarrow{\partial}
 D_J\longrightarrow0,\qquad
 D_J=\frac{J^{m-1}\cap\im\delta}{\im\delta_J}.
\end{equation}
For two subalgebras,
$D_J=(J\cap(A_1+A_2))/((J\cap A_1)+(J\cap A_2))$.
For an arbitrary family the same formula holds with products in place of direct sums and a chosen base index.
\end{theorem}
\begin{proof}
Lift a common image to $(a_i)$ and send it to $[\delta(a_i)]$. Its components lie in $J$. Changing lifts changes this by $\im\delta_J$, and every numerator element occurs. The class is zero precisely when subtracting elements of $A_i\cap J$ makes the tuple diagonal, precisely the condition that the common image lifts to $C$. This proves exactness, including the product version.
\end{proof}

\begin{corollary}\label{allowed-correction}
If an allowed subalgebra $U$ is contained in $K$, then
\[
 0\longrightarrow\frac{q(C)}{q(C)\cap U}\longrightarrow\frac K U
 \longrightarrow\frac{D_J}{\partial(U)}\longrightarrow0
\]
is exact. A sufficient condition for $D_J=0$ is a linear map $N:D\to D$ with $\ker N=J$, $qN=q$, and $N(A_i)\subseteq A_i$ for every $i$.
\end{corollary}
\begin{proof}
Pass \eqref{descent-exact} to the quotient by $U$. For the sufficient condition, lifts of a common image have differences in $J$, so their $N$-images coincide and give a common lift.
\end{proof}
If $t\in Z(D)$ acts injectively on $D$, and $A\cap tD=tA$, $B\cap tD=tB$, then the two-chart defect for $J=tD$ is the $t$-torsion of $D/(A+B)$. Indeed multiplication by $t$ identifies its numerator and denominator with $\{d:td\in A+B\}$ and $A+B$, respectively.

\begin{example}
In $D=\C\alg{x^{\pm1},y^{\pm1}}$ take
$A=\C\alg{x^{\pm1},y}$ and $B=\C\alg{x,y^{\pm1}}$. Reduced-word bases give $A\cap B=\C\alg{x,y}$. Under $y=x$ both images are $\C[x^{\pm1}]$, and the defect is
\[
 \C[x^{\pm1}]/\C[x]=\bigoplus_{n\geq1}\C x^{-n}.
\]
The numerator representative of $x^{-n}$ in the two-chart formula is $x^{-n}-y^{-n}$.
\end{example}

\begin{theorem}[Tensor-layer intersections]\label{tensor-layer}
Let $M$ be a flat right $S$-module, $N$ a flat left $S$-module, $A\subseteq M$ a right submodule, and $B\subseteq N$ a left submodule. Then
\begin{equation}\label{tor-intersection}
 \frac{(A\otimes_S N)\cap(M\otimes_S B)}
 {\im(A\otimes_S B\to M\otimes_S N)}
 \simeq\Tor^S_1(M/A,N/B).
\end{equation}
The two modules in the numerator are identified with their images in $M\otimes_S N$.
\end{theorem}
\begin{proof}
Flatness embeds those two modules. The kernel of $A\otimes N\to A\otimes(N/B)$ is the image of $A\otimes B$. Consequently the quotient on the left is the kernel of
$A\otimes(N/B)\to M\otimes(N/B)$.
The long exact sequence for $0\to A\to M\to M/A\to0$, using flatness of $M$, identifies that kernel with the indicated Tor group.
\end{proof}

\begin{proposition}[Free-group boundary complex]\label{boundary-complex}
Let $S=k\Fgroup_d$ with free generators $z_i$. For a right $S$-module $Q$ and left $S$-module $V$,
\begin{equation}\label{boundary-formula}
 \Tor^S_1(Q,V)=\ker\left[\bigoplus_{i=1}^d(Q\otimes_k V)\longrightarrow Q\otimes_k V\right],
 \quad q\otimes v\longmapsto qz_i\otimes v-q\otimes z_i v,
\end{equation}
and $\Tor^S_n(Q,V)=0$ for $n\geq2$.
\end{proposition}
\begin{proof}
There is a free right-module resolution
\[
 0\to\bigoplus_i(Q\otimes_k S)\xrightarrow{d_1}Q\otimes_k S\to Q\to0,
 \qquad d_1(q\otimes s)=qz_i\otimes s-q\otimes z_i s.
\]
The last map is multiplication. To check exactness, change vertex coordinates by $q\otimes g\mapsto qg\otimes[g]$ and edge coordinates by $q\otimes g\mapsto qz_i g\otimes[i,g]$. The resulting complex is the augmented finite-support chain complex of the Cayley tree with vector-space coefficients in $Q$. A leaf proves injectivity on edge chains; paths from a fixed vertex express every augmentation-zero vertex chain as a boundary. Tensoring with $V$ proves \eqref{boundary-formula} and higher vanishing.
\end{proof}
For example $S=k[z^{\pm1}]$, $M=N=S$, and $A=B=(z-1)S$ give a nonzero correction
$(z-1)S/(z-1)^2S\simeq k$.
For a surface presentation, the quotient formula uses the relation ideal. The tensor formula requires a filtration whose tensor layers satisfy the stated flatness hypotheses.

\section{Noncommutative polygons}\label{sec-polygons}
For a polygon with cyclic vertices $[n]$, the Berenstein--Retakh algebra has invertible directed chord generators $x_{ij}$ and relations
\begin{align}\label{polygon-relations}
 x_{ij}x_{kj}^{-1}x_{ki}&=x_{ik}x_{jk}^{-1}x_{ji},\\
 x_{jl}&=x_{jk}x_{ik}^{-1}x_{il}+x_{ji}x_{ki}^{-1}x_{kl},\notag
\end{align}
where the first has distinct vertices and the second cyclically ordered $i,j,k,l$. The chart $\A_\Delta$ is generated by directed triangulation edges and their inverses; $\U_n$ is generated by all chords and inverses of the directed boundary edges. We use the chart and Laurent theorems of \cite[Theorems 2.8 and 2.10]{BR}: the chart is the group algebra of its triangle presentation, it contains all chords, and the full algebra is its universal localization at all chords. In particular $\U_n\subseteq\A_\Delta$.

\begin{lemma}[Boundary normalization]\label{gauge}
Suppose the directed boundary edges to be normalized form a forest in the bipartite graph with vertices $L_i,R_i$ and edge $L_iR_j$ for $x_{ij}$. Setting those boundary coordinates to one gives compatible decompositions of the full algebra, every chart, and the allowed algebra as the free product of the normalized algebra with the free group algebra on the forest edges.
\end{lemma}
\begin{proof}
Root each component and solve $w_{ij}=g_i h_j$ successively along its edges, fixing a root factor to one. The factors are words in the independent boundary units. Substitute $x_{ij}=g_iq_{ij}h_j$. In a triangle relation the internal factors cancel and both sides have outer factors $g_i,h_i$. In exchange both summands have the same outer factors $g_j,h_l$. Thus the remaining relations are exactly the normalized ones and impose no relation on the forest letters. The inverse substitution reconstructs every original generator. Apply the same substitution to the abstract triangle presentation to obtain the compatible chart decompositions. Restriction of its embedded chart map proves injectivity on the normalized factor. Boundary units also recover the allowed generators. In the surface case a corner angle transforms by $h_z^{-1}\widehat\Theta_z g_z^{-1}$, so the statement applies with total angles as well.
\end{proof}

\begin{lemma}\label{freeproduct-intersections}
For unital subalgebras $R_\alpha\subseteq D$ and a unital $k$-algebra $H$,
\[
 \bigcap_\alpha(H*R_\alpha)=H*\bigl(\bigcap_\alpha R_\alpha\bigr)
 \quad\text{inside }H*D.
\]
\end{lemma}
\begin{proof}
Choose $H=k\oplus H^\circ$, $D=k\oplus D^\circ$. In each alternating tensor normal-form summand, membership replaces every $D^\circ$ slot by $R_\alpha\cap D^\circ$. Tensor products over a field commute with intersections in one slot, and
$(V'\otimes W)\cap(V\otimes W')=V'\otimes W'$. Apply this successively to the finitely many slots. Arbitrary intersections are allowed because any given tensor has finite support; equivalently coefficient slice maps test membership in each slot.
\end{proof}

For $n\geq4$ put $m=n-3$. Normalize both orientations of the boundary path $1,2,\ldots,n$, leaving $x_{1n},x_{n1}$ unnormalized. Its $2n-2$ directed edges form a bipartite forest with two components. In the fan at vertex $1$, set
\[
 u_0=1,\quad u_i=x_{1,i+2}\ (1\leq i\leq m+1),\quad
 t=u_1,\quad d_i=x_{i+1,i+3}.
\]
Define right and left continuants by
\begin{equation}\label{polygon-continuants}
 v_0=\rho_0=1,\quad v_1=\rho_1=t,\quad
 v_{i+1}=v_i d_i-v_{i-1},\quad \rho_{i+1}=d_i\rho_i-\rho_{i-1},\quad S=v_{m+1}.
\end{equation}

\begin{lemma}[Tridiagonal continuants]\label{tridiagonal}
Over an arbitrary unital ring, let $T_j$ be tridiagonal of size $j$ with diagonal $t,d_1,\ldots,d_{j-1}$ and adjacent entries one. Then $T_j$ is invertible if and only if $v_j$ is a unit, and
\[
 P\alg{T_j^{-1}}=P[v_j^{-1}].
\]
For $j=m+1$, $\rho_{m+1}$ is also a unit in $P[S^{-1}]$, and
\begin{equation}\label{opposite-continuant}
 \rho_{m+1}^{-1}=K_m S^{-1}v_m-K_{m-1},
\end{equation}
where $K_i=K(d_1,\ldots,d_i)$, $K_0=1$, and the right continuant satisfies
\[K(a_1,\ldots,a_i)=K(a_1,\ldots,a_{i-1})a_i-K(a_1,\ldots,a_{i-2}).\]
\end{lemma}
\begin{proof}
For $j=m+1$, the row $(1,-v_1,v_2,\ldots,(-1)^m v_m)$ times $T_j$ is $(0,\ldots,0,(-1)^m S)$. Solving the last $m$ equations of $T_jx=b$ backwards expresses all entries in terms of $x_m$ and gives
\[
 Sx_m=\sum_{i=0}^m(-1)^{m+i}v_i b_i.
\]
If $S$ is a unit this gives a unique solution with inverse entries in $P[S^{-1}]$, and $(T_j^{-1})_{m0}=(-1)^mS^{-1}$. Conversely invertibility of $T_j$ makes left multiplication by $S$ surjective, by $b=e_0$, and injective, by backward substitution of a vector with $Sx_m=0$. Hence a right inverse of $S$ is also a left inverse. This proves the equivalence over every unital ring and the inverse-algebra equality.

For $b=e_m$, one has $x_m=S^{-1}v_m$ and backward substitution gives
$x_0=(-1)^m(K_mS^{-1}v_m-K_{m-1})$. Solving the first $m$ homogeneous equations forwards gives $x_i=(-1)^i\rho_i x_0$, so the last equation says $\rho_{m+1}(K_mS^{-1}v_m-K_{m-1})=1$. If $\rho_{m+1}z=0$, the vector $((-1)^i\rho_i z)_i$ lies in the kernel of $T_j$, so $z=0$. The right inverse is thus two-sided. The argument applies at every truncation $j$. The backward continuant formula follows equally by its first-end recurrence, obtained by grouping path matchings by their first vertex.
\end{proof}

\begin{proposition}\label{polygon-model}
In the normalized fan chart,
\[
 \widehat\A_F=\C\alg{u_1^{\pm1},\ldots,u_{m+1}^{\pm1}},\qquad
 u_i=v_i,\qquad x_{n1}=\rho_{m+1},\qquad
 \widehat\U_n=P[S^{-1}],\quad P=\C\alg{t,d_1,\ldots,d_m}.
\]
This free group chart embeds faithfully in the free skew field of $P$.
\end{proposition}
\begin{proof}
Put $b_j=x_{j+2,1}$, with $b_0=1$. The fan triangle equations are
$b_{j+1}=u_j^{-1}u_{j+1}b_j$.
They eliminate all reverse edges, impose no relations on the $u_j$, and prove the group-chart presentation. Exchange and the triangle relation give
$u_i d_i=u_{i-1}+u_{i+1}$ and $x_{i+3,i+1}=d_i$.
Then $b_1=t$ and $b_{i+1}=d_i b_i-b_{i-1}$, so $b_i=\rho_i$.
For chords not incident with $1$ the same quadrilateral identities give
\[
 x_{a,j+1}=x_{aj}d_{j-2}-x_{a,j-1},\qquad
 x_{j+1,a}=d_{j-2}x_{ja}-x_{j-1,a},
\]
starting with consecutive values one and two-step values $d_{a-1}$. Every chord is polynomial in $P$. Only the two opposite boundary inverses remain; Lemma~\ref{tridiagonal} supplies the second from $S^{-1}$.

The mutually inverse rational substitutions are $u_i=v_i$, $t=u_1$, $d_i=u_i^{-1}(u_{i-1}+u_{i+1})$. At
\begin{equation}\label{polygon-center}
 t=1,\qquad d_i=2
\end{equation}
all $v_i=1$. The linear terms $V_i=v_i-1$ satisfy
$V_1=Y_0$, $V_{i+1}=2V_i+Y_i-V_{i-1}\pmod{(Y)^2}$.
This triangular linear change is invertible, so substitution is an isomorphism on completed free power-series rings, by solving in each homogeneous degree. The Magnus embedding of the free group chart followed by this substitution is faithful. A zero rational function has zero scalar-center series, proving the asserted embedding.
\end{proof}
At \eqref{polygon-center}, all chords incident with $1$ equal one; the other chords have value $|b-a|$ for endpoints $a,b\geq2$. All are nonzero. Thus the universal localization of the fan chart maps to this skew field; faithfulness is needed only on the fan chart when pulling an intersection equality back.

\begin{theorem}[Polygon upper algebra]\label{polygon-theorem}
For $k=\Q$ or $\C$ and $n\geq4$, let $F(i)$ be the flip of $(1,i+2)$ in the fan. Then
\begin{equation}\label{polygon-upper}
 \U_n=\A_F\cap\bigcap_{i=1}^{n-3}\A_{F(i)}=\bigcap_\Delta\A_\Delta.
\end{equation}
For $n=2,3$ the last equality also holds.
\end{theorem}
\begin{proof}
Normalize as above and work over $\C$. At \eqref{polygon-center}, the inverse $C_j$ of $T_j$ has entries
\begin{equation}\label{tridiagonal-inverse-center}
 (C_j)_{ab}=(-1)^{a+b}(j-\max(a,b)),\qquad 0\leq a,b<j.
\end{equation}
This follows by multiplying by the tridiagonal matrix; in particular every entry is nonzero. The coefficient matrices $C_jE_{aa}$ of its monic normalization generate $\Mat_j(\C)$, since
\[
 (C_jE_{aa})(C_jE_{bb})=(C_j)_{ab}C_jE_{ab}.
\]
Thus $v_j^{-1}$ has one class $V_j$, of dimension $j$, and the classes are distinct. The fan support is $\{V_1,\ldots,V_{m+1}\}$ and $P[S^{-1}]$ has support $\{V_{m+1}\}$.

The $i$th flip has directed generators $u_j^{\pm1}$ for $j\neq i$ and $d_i^{\pm1}$. For precision, the reverse edges before the omitted position are eliminated as above; the triangle spanning the flipped quadrilateral gives
\[
 b_{i+1}=d_i u_{i-1}^{-1}u_{i+1}d_i^{-1}b_{i-1},
\]
and the same triangle recursion recovers the subsequent reverse fan edges. For $i=1$ use $u_0=b_0=1$. Hence there are no further inverse generators. The scalar inverse class of $d_i$ differs from $V_i$, by dimension if $i>1$, and by its coefficient variable if $i=1$. That flip removes $V_i$. Every chart contains $P[S^{-1}]$, so the common support of these $m+1=n-2$ charts is exactly $\{V_{m+1}\}$. Theorem~\ref{support-intersections} proves the normalized equality in the skew field. Faithfulness on the fan pulls it back. Lemmas~\ref{gauge} and~\ref{freeproduct-intersections} restore the boundary free factor.

The full chart intersection is contained in this finite one and contains $\U_n$ by the Laurent theorem. This proves \eqref{polygon-upper}. For $n=2,3$ every chord is a boundary edge. All data are defined over $\Q$; scalar extension to $\C$ commutes with finite intersections of vector subspaces and is faithful. The finite equality therefore descends to $\Q$, and the same two inclusions give the full intersection there.
\end{proof}
This is the assertion of \cite[Conjecture 2.20]{BR}, together with its finite fan-and-flips form.

\section{Unpunctured M\"obius strips}\label{sec-mobius}
Let $M_n$ be a M\"obius strip with $n$ marked boundary points and no punctures. A directed curve is a relative homotopy class with endpoints in the marked set, with the trivial class assigned one. Reversal is denoted by a bar. The unpunctured surface algebra $\A_{M_n}$ has invertible nontrivial curve generators and the triangle and exchange relations \eqref{polygon-relations} on developed polygons, projected to the surface. These are the unpunctured conventions of \cite[Sections 3.3--3.4]{BR}. An ordinary triangulation consists of simple based arcs with disjoint interiors, including both orientations and the boundary segments. No unbased closed core curve is included.

For a marked point $z$, the angle subtended by $a,b$ in a developed triangle is
\[
 \Theta_z(a,b)=x_{az}^{-1}x_{ab}x_{zb}^{-1}.
\]
The total angle $T_z$ is the sum over its incident triangle sectors. Let $\U_{M_n}$ be generated by all curve elements, inverses of directed boundary segments, and the total angles. We use the surface chart, Laurent, and angle theorems \cite[Theorems 3.30, 3.36, 3.40]{BR}: the triangle group chart embeds, every curve belongs to every chart, the full algebra is its universal curve-inverting localization, and total angles are chart-independent. These facts give $\U_{M_n}\subseteq\A_\Delta$ without assuming injectivity of the full localization into a skew field.

\subsection{The covering triangulation and free coordinates}
Use the strip $\R\times[0,1]$ with deck map $g(x,y)=(x+1,1-y)$. Put
\[
 p_0=(0,0),\quad p_1=(1,1),\quad p_2=(2,0),\quad
 r_j=(2j/n,0)\ (0\leq j\leq n-1).
\]
The triangular fundamental region with boundary vertices
$(p_0,r_1,\ldots,r_{n-1},p_2,p_1)$
has its two slanted sides identified by $g$. Its fan at $p_1$, together with the seam $p_0p_1$, projects to a triangulation $T$ with $n$ interior arc orbits and $n$ triangles. Deck translates tile the strip with disjoint triangle interiors.

Normalize both orientations of $r_{j-1}r_j$, $1\leq j\leq n-1$, to one. These directed edges form a bipartite forest with two components, so Lemma~\ref{gauge} applies. Put
\[
 v=x_{p_0p_1},\quad u=\bar v,\quad w_j=x_{r_jp_1},\quad t=w_1,\quad s=\bar t,\quad b=x_{r_{n-1}p_2},\quad m=n-2.
\]
The normalized fan chart is free on the group basis
\begin{equation}\label{mobius-fan-basis}
 s,v,w_1,\ldots,w_{n-1},b.
\end{equation}
Indeed the first triangle gives $u=svt^{-1}$; subsequent triangles eliminate reverse fan edges by
$\bar w_j=\bar w_{j-1}w_jw_{j-1}^{-1}$,
with $w_0=v$, $\bar w_0=u$. The terminal triangle eliminates the opposite boundary value $c=\bar b$ by
$c=uw_{n-1}^{-1}bv^{-1}\bar w_{n-1}$.
These eliminate all other directed edges and impose no remaining relations on \eqref{mobius-fan-basis}.

Let
\[
 d_j=x_{r_{j-1}r_{j+1}}=\bar d_j,\qquad
 w_{j+1}=d_jw_j-w_{j-1},\qquad 1\leq j\leq m.
\]
The equal two-step orientations and the recurrence follow from the triangles on consecutive boundary marks and their adjacent quadrilaterals, exactly as in \eqref{polygon-relations}. Let $q_1=gr_1$, $q=x_{r_1q_1}$, $p=\bar q$, and $r=x_{p_0q_1}$.

\begin{lemma}[Seed and boundary identities]\label{mobius-seed}
In the actual normalized fan chart,
\begin{equation}\label{mobius-seed-eq}
 u=svt^{-1},\quad tp=qs,\quad t\bar r=1+qu,\quad
 r=vt^{-1}q+s^{-1}.
\end{equation}
Consequently
\begin{equation}\label{mobius-H}
 H=\begin{pmatrix}-t&-q\\v&r\end{pmatrix},\quad
 H^{-1}=\begin{pmatrix}-\bar r&-p\\u&s\end{pmatrix},\quad
 H^*JH=-J,\qquad J=\begin{pmatrix}0&1\\-1&0\end{pmatrix},
\end{equation}
where $*$ is transpose followed by reversal. Define
\begin{equation}\label{frames}
 F_0=e_1,\quad F_1=e_2,\quad F_{j+1}=F_jd_j-F_{j-1},\quad
 F_j=\binom{-A_j}{B_j},\quad G_j=(F_{j-1},F_j).
\end{equation}
Then $G_j^{\pm1}$ are polynomial, $G_j^*JG_j=J$, and
\begin{equation}\label{mobius-boundary}
 w_j=F_j^*JHe_1,\qquad b=F_{n-1}^*JH^2e_1.
\end{equation}
\end{lemma}
\begin{proof}
The triangles $(p_0,r_1,p_1)$ and $(r_1,p_1,q_1)$ give the first two identities in \eqref{mobius-seed-eq}; the cyclic quadrilateral $(r_1,q_1,p_1,p_0)$ gives the third. Reversal gives
$r=(1+vp)s^{-1}=s^{-1}+vt^{-1}q$.
Substitution proves both matrix inverse products in \eqref{mobius-H}, and $H^{-1}=JH^*J$ gives the anti-unitary identity. The transition
$G_{j+1}=G_j\left(\begin{smallmatrix}0&-1\\1&d_j\end{smallmatrix}\right)$
proves the frame assertions. The first pairing in \eqref{mobius-boundary} is the left recurrence with initial values $v,t$.

For the boundary identity, use actual curves $\theta_j=x_{r_j,q_1}$. The quadrilateral $(r_{j-1},r_j,r_{j+1},q_1)$ gives the same left recurrence, with $\theta_0=r$, $\theta_1=q$. Thus
$\theta_j=\bar B_jq-\bar A_jr$.
In $(r_{n-1},p_2,q_1,p_1)$ the sides $x_{p_2q_1}=s$, $x_{p_2p_1}=u$, $x_{q_1p_2}=t$, and $x_{q_1p_1}=1$ give
$\theta_{n-1}s^{-1}u=w_{n-1}+bt^{-1}$.
Use $s^{-1}u=vt^{-1}$ and multiply by $t$ to get $\theta_{n-1}v=b+w_{n-1}t$. With $A=\bar A_{n-1}$, $B=\bar B_{n-1}$ this reads
\begin{equation}\label{b-polynomial}
 b=-A(rv-vt)+B(qv-t^2),
\end{equation}
which is exactly $F_{n-1}^*JH^2e_1$. This derivation precedes the general curve-pairing construction.
\end{proof}

\begin{proposition}[Free polynomial model]\label{mobius-free-model}
The variables $t,q,v,r,d_1,\ldots,d_m$ are free polynomial coordinates, and their free skew field is the skew field of the normalized fan chart. Their common scalar center corresponding to all generators in \eqref{mobius-fan-basis} equal to one is
\begin{equation}\label{mobius-alpha}
 t=v=1,\qquad d_j=2,\qquad q=n,\qquad r=n+1.
\end{equation}
\end{proposition}
\begin{proof}
From the group coordinates set
\[
 t=w_1,\quad d_j=(w_{j-1}+w_{j+1})w_j^{-1},\quad
 q=tw_{n-1}^{-1}(b+w_{n-1}t+As^{-1}v)v^{-1},\quad
 r=vt^{-1}q+s^{-1}.
\]
Conversely recover the $w_j$ by their recurrence, $s=(r-vt^{-1}q)^{-1}$, and $b$ by \eqref{b-polynomial}. These formulas are inverse, since
$b=w_{n-1}t^{-1}qv-As^{-1}v-w_{n-1}t$.
They are regular at the specified centers. Their inverse substitutions give an isomorphism of completed free power-series algebras. Faithfulness of the Magnus expansion of the group chart proves polynomial freeness, and the inverse rational substitutions give the equality of free skew fields.
\end{proof}

\subsection{All winding curves and total angles}
Put
\begin{equation}\label{mobius-rings}
 P_M=\C\alg{t,q,v,r,d_1,\ldots,d_m},\qquad
 R_0=P_M\alg{H^{-1}},\qquad R_M=R_0[b^{-1}].
\end{equation}
The reversal on the chart extends to its skew field. It fixes $d_j$ and preserves $R_0$, by \eqref{mobius-H}.

\begin{lemma}[Signed pairings]\label{signed-pairings}
Over a skew field with involution, let $E_i$ be isotropic columns for $J$ with $E_i^*JE_j\neq0$ for distinct vertices. Order the vertices along the boundary of a developed polygon and put
\[
 X_{ij}=\sigma_{ij}E_i^*JE_j,\qquad \sigma_{ij}=1\ (i<j),\quad\sigma_{ij}=-1\ (i>j).
\]
Then $X_{ji}=\overline{X_{ij}}$, and the directed triangle and exchange identities hold. For columns $H^kF_j$ assigned to $g^kr_j$, these values are invariant under deck translation.
\end{lemma}
\begin{proof}
The reversal identity follows from $J^*=-J$. For $i<j<k$, normalize the transverse isotropic columns $E_i,E_k$ to $e_1,e_2$. This is possible because their pairing is nonzero: rescale the second column to make their pairing one, after which their frame is symplectic. Right frame rescalings transform $X_{ij}$ by $G_i^*X_{ij}G_j$, preserving the directed relations. Write $E_j=(A,B)^T$. Isotropy says $A^*B=B^*A$, and the triangle relation becomes
$BA^{-1}=(A^*)^{-1}B^*$,
which is that identity multiplied by inverses. Other orders follow by multiplication and inversion.
For exchange, expand a fourth column in the frame $(E_i,E_k)$. Writing $P_{ab}=E_a^*JE_b$ gives
\[
 P_{jl}=P_{ji}P_{ki}^{-1}P_{kl}+P_{jk}P_{ik}^{-1}P_{il}.
\]
The boundary-order signs convert this to \eqref{polygon-relations}; cyclic rotation and reversal cover all orientations. Finally $g$ reverses the boundary order of the strip and $H$ is anti-unitary, so its sign and the pairing sign cancel.
\end{proof}

\begin{lemma}\label{mobius-nonzero}
At \eqref{mobius-alpha} all pairings of distinct lifted marked vertices are nonzero.
\end{lemma}
\begin{proof}
There $F_j=(1-j,j)^T$ and
$H_\alpha=\left(\begin{smallmatrix}-1&-n\\1&n+1\end{smallmatrix}\right)$,
with trace $n$ and determinant $-1$. For $0\leq i,j\leq n-1$, put $a_k=F_i^TJH_\alpha^kF_j$. Direct multiplication and Cayley--Hamilton give
\[
 a_0=j-i,\qquad a_1=1+j(n-i)\geq1,\qquad a_{k+2}=na_{k+1}+a_k.
\]
Since $a_0\geq-(n-1)$, $a_2\geq1$, and all $a_k$ with $k\geq1$ are positive. Negative powers follow from
$a_{-k}(i,j)=(-1)^{k+1}a_k(j,i)$.
At $k=0$ only an identical lift gives zero. Apply deck translation to reduce every pairing to this calculation.
\end{proof}

\begin{proposition}[Surface-coordinate identification]\label{mobius-coordinate}
Every nontrivial curve equals its signed pairing of endpoint columns $H^kF_j$ in the fan skew field. In particular every curve belongs to $R_0$, and
\begin{equation}\label{mobius-allowed}
 \widehat\U_{M_n}=R_M.
\end{equation}
\end{proposition}
\begin{proof}
The identity $G_j^*JG_j=J$ gives the normalized boundary pairings. Equations \eqref{mobius-boundary} and \eqref{b-polynomial} give the seam, fan, and terminal boundary values. Lemma~\ref{signed-pairings} gives agreement under reversal and deck translation. The pairings are nonzero by Lemma~\ref{mobius-nonzero}.

Develop the finite chain of triangles crossed by a nontrivial relative curve into its canonical polygon \cite[Theorem 3.21 and Lemma 3.23]{BR}. The distinct lifted endpoints determine the reduced class in this strip. The pairings give the values of the triangulation edges in the developed polygon and satisfy its relations by Lemma~\ref{signed-pairings}. The polygon Laurent expression therefore identifies the diagonal pairing with the curve \cite[Theorem 3.30]{BR}. The pairing is polynomial in $H^{\pm1}$ and the $d_j$, so it belongs to $R_0$. The generators of $R_0$ are curve values in the fan chart by \eqref{mobius-H}, and the faithful chart map pulls back these identities. Since each curve has a nonzero image in the skew field, the map extends to the universal surface localization.

It remains to check boundary inverses and angles. Put $G=G_{n-1}$, $M=G^{-1}H^2$, and $c=\bar b$. Then
\[
 M_{11}=-b,\qquad (M^{-1})_{22}=-c.
\]
The first equality is \eqref{mobius-boundary}; the second follows from $M^*JM=J$ and reversal. The Schur complement gives
\begin{equation}\label{reverse-boundary-inverse}
 c^{-1}=-M_{22}-M_{21}b^{-1}M_{12}\in R_M.
\end{equation}
All other normalized boundary inverses are one.

For total angles, exchange on the cyclic quadruple $(z,a,b,c)$ gives
\[
 x_{ac}=x_{ab}x_{zb}^{-1}x_{zc}+x_{az}x_{bz}^{-1}x_{bc}.
\]
Multiplication by $x_{az}^{-1}$ and $x_{zc}^{-1}$ yields
$\Theta_z(a,c)=\Theta_z(a,b)+\Theta_z(b,c)$.
Develop the finite star of a boundary lift $z$ and telescope. If $\ell,z,r'$ are consecutive boundary lifts,
\begin{equation}\label{angle-telescope}
 T_z=x_{\ell z}^{-1}x_{\ell r'}x_{zr'}^{-1}.
\end{equation}
The outer inverses are boundary inverses, and the middle curve belongs to $R_0$. Thus all allowed generators belong to $R_M$. Conversely every generator of $P_M$, every entry of $H^{-1}$, and $b^{-1}$ is an allowed curve or boundary inverse, proving \eqref{mobius-allowed}.
\end{proof}

\subsection{Ordinary flips and the finite intersection}
For any two lifted endpoints choose invertible symplectic or anti-symplectic frames whose selected columns are their columns. They can be obtained from $H^kG_j$ and constant basis permutations. Their pairing matrix is invertible over $R_0$, and the curve and its reversal are opposite inverse entries up to sign. Lemma~\ref{opposite} therefore shows that inverting one directed curve inverts its reversal over $R_0$. Consequently the image of the normalized fan chart is
\begin{equation}\label{mobius-fan-localization}
 R_M[v^{-1},w_1^{-1},\ldots,w_{n-1}^{-1}].
\end{equation}
The same observation describes a flipped chart by retaining the other inverses and adjoining one replacement inverse. The ring $R_M$ is a finite pencil inverse algebra: take a block sum of $H$ and a polynomial Schur linearization of $b$, normalized at \eqref{mobius-alpha}. Write $C$ for its support.

\begin{lemma}\label{mobius-tridiagonal}
For $1\leq j\leq n-1$, $P_M[w_j^{-1}]$ is the inverse algebra of the affine tridiagonal matrix
\begin{equation}\label{mobius-tridiagonal-eq}
 T_j=\begin{pmatrix}
 t&1&0&\cdots&0\\
 v&d_1&1&\ddots&\vdots\\
 0&1&d_2&\ddots&0\\
 \vdots&\ddots&\ddots&\ddots&1\\
 0&\cdots&0&1&d_{j-1}
 \end{pmatrix},\qquad T_1=(t).
\end{equation}
Its monic normalization is irreducible of dimension $j$. Denote its class by $W_j$. The $W_j$ are pairwise distinct, and the coefficient of $d_{j-1}$ is nonzero for $j\geq2$. The scalar class $V$ of $v$ differs from every $W_j$.
\end{lemma}
\begin{proof}
Solve the first $j-1$ equations forwards; all unknowns are determined by the first, and the coefficient in the last equation is $(-1)^{j-1}w_j$. As in Lemma~\ref{tridiagonal}, this proves two-sided invertibility equivalence over an arbitrary ring and gives $(T_j^{-1})_{0,j-1}=(-1)^{j-1}w_j^{-1}$. At \eqref{mobius-alpha} the tridiagonal matrix is the one of \eqref{tridiagonal-inverse-center}. Its diagonal-variable coefficient matrices generate $\Mat_j(\C)$ by the same rank-one product calculation. The additional coefficient of $v$ does not alter this. Dimension distinguishes the $W_j$; the last diagonal coefficient is nonzero. Finally $V$ differs from $W_1$, the class of $t$, by its coefficient variable, and from the others by dimension.
\end{proof}

\begin{lemma}[Covering quadrilaterals]\label{mobius-flips}
For $n\geq2$ every interior arc orbit of $T$ has an ordinary flip. The replacement coordinates, in suitable orientations, are
\begin{equation}\label{flip-list}
 \begin{array}{c|c|c}
 \text{removed arc}&\text{replacement lift}&\text{coordinate}\\\hline
 r_jp_1\ (1\leq j\leq n-2)&r_{j-1}r_{j+1}&d_j\\
 p_0p_1&r_{n-1},gr_1&\theta=F_{n-1}^*JHe_2\\
 r_{n-1}p_1&r_{n-2}p_2&\eta=F_{n-2}^*JH^2e_1
 \end{array}
\end{equation}
All replacement coordinates equal two at \eqref{mobius-alpha}.
\end{lemma}
\begin{proof}
The lifted fan triangles are $(r_{j-1},r_j,p_1)$, $1\leq j\leq n-1$, and $(r_{n-1},p_2,p_1)$. Two adjacent ones give the first quadrilateral, and the last two give the final quadrilateral. Across the seam develop the last triangle by $g^{-1}$; it and the first triangle give opposite vertices $r_1,g^{-1}r_{n-1}$. Translate by $g$ and choose the reverse orientation to obtain the middle lift in \eqref{flip-list}.

The displayed strip triangulation gives the ordinary-arc assertion directly. Each lifted triangle has at most one side in the orbit of the selected arc when $n\geq2$. Distinct deck translates of its two-triangle quadrilateral therefore have disjoint triangle interiors. The new diagonal, drawn inside the two open triangles and crossing only their common side, has translates with disjoint interiors and misses every retained edge. It projects to a simple based arc, with a crossing of the removed arc, and hence is its ordinary flip. A same-boundary replacement is drawn through the open quadrilateral, avoiding intervening marked boundary points. This argument also applies for $n=2$, where the two quotient triangles share two distinct arc orbits: their lifted quadrilaterals still have disjoint interiors. At $n=1$ the one triangle has two sides in the seam orbit, and this argument does not apply. The coordinate formulas follow from Proposition~\ref{mobius-coordinate}. Evaluating their pairings by the recurrence in Lemma~\ref{mobius-nonzero} gives two.
\end{proof}

\begin{figure}[ht]
\centering
\begin{tikzpicture}[x=2.1cm,y=2cm, every node/.style={font=\small}]
\fill[black!6] (0,0)--(1,0)--(1,1)--(0,1)--cycle;
\draw (-.1,0)--(2.1,0) (-.1,1)--(2.1,1);
\draw (0,1)--(0,0) (1,0)--(1,1)--(2,0);
\draw[very thick] (0,0)--(1,1);
\draw[dashed,thick] (1,0)--(0,1);
\foreach \p in {(0,0),(1,0),(2,0),(0,1),(1,1)} \fill \p circle (1.1pt);
\node[below] at (0,0) {$p_0$};
\node[below] at (1,0) {$r_1$};
\node[below] at (2,0) {$p_2$};
\node[above] at (0,1) {$g^{-1}r_1$};
\node[above] at (1,1) {$p_1$};
\node at (1,-.48) {Seam flip};
\begin{scope}[xshift=6.2cm]
\fill[black!6] (0,0)--(2,0)--(1,1)--cycle;
\draw (-.1,0)--(2.1,0) (-.1,1)--(2.1,1);
\draw (0,0)--(1,1)--(2,0);
\draw[very thick] (1,0)--(1,1);
\draw[dashed,thick] (0,0).. controls (.6,.55) and (1.4,.55)..(2,0);
\foreach \p in {(0,0),(1,0),(2,0),(1,1)} \fill \p circle (1.1pt);
\node[below] at (0,0) {$p_0$};
\node[below] at (1,0) {$r_1$};
\node[below] at (2,0) {$p_2$};
\node[above] at (1,1) {$p_1$};
\node at (1,-.48) {Last fan flip};
\end{scope}
\end{tikzpicture}
\caption{Developed quadrilaterals for $n=2$. Heavy solid arcs are removed and dashed arcs replace them. The shaded regions consist of the two incident lifted triangles. The right replacement avoids the intervening boundary mark.}
\label{fig-mobius-flips}
\end{figure}
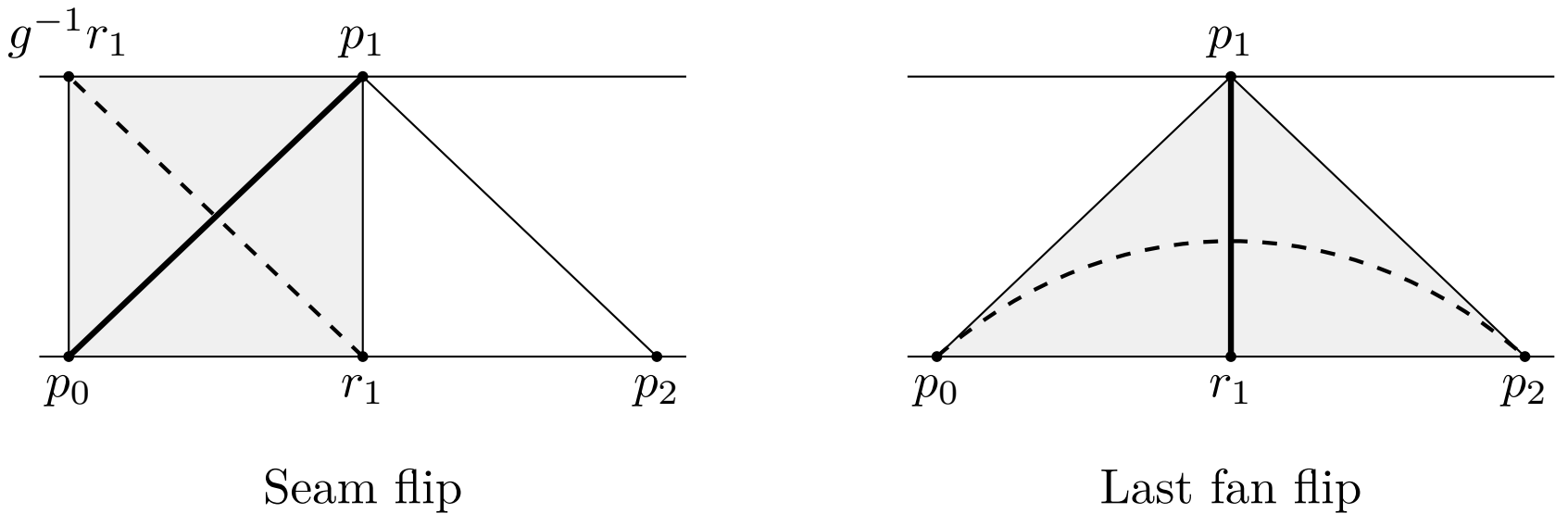

\begin{theorem}[M\"obius upper algebra]\label{mobius-upper}
For an unpunctured M\"obius strip with $n\geq2$ boundary marks, over $\Q$ or $\C$,
\begin{equation}\label{mobius-upper-eq}
 \U_{M_n}=\A_T\cap\bigcap_{e\in T\text{ interior}}\A_{\mu_eT}
          =\bigcap_\Delta\A_\Delta.
\end{equation}
The finite intersection uses $n+1$ ordinary charts.
\end{theorem}
\begin{proof}
Work first in the normalized free coordinates over $\C$. By Lemma~\ref{mobius-tridiagonal} the fan support is
$C\cup\{V,W_1,\ldots,W_{n-1}\}$.
For $1\leq j\leq n-2$ its flip replaces $w_j$ by $d_j$, whose scalar inverse class differs from $W_j$ by dimension, or by variable when $j=1$. All retained classes are distinct, so this removes $W_j$ outside $C$.

The seam replacement is
$\theta=\bar B_{n-1}q-\bar A_{n-1}r$;
it omits $v$. A polynomial linearization can omit that variable as well, so none of its constituents can be the scalar class $V$, whose $v$ coefficient is nonzero. For $n\geq3$, the final replacement $\eta$ omits $d_{n-2}$, whereas $W_{n-1}$ has a nonzero coefficient of that variable. It therefore removes $W_{n-1}$ outside $C$. The omitted-variable statement is invariant under constant similarity, so it indeed excludes the indicated simple classes.

For $n=2$, the final replacement is $\eta=rv-vt$. Its inverse support cannot contain the scalar class of $t$: otherwise the constituent decomposition would imply $Z(t)\subseteq Z(L_\eta)$, and since $\eta$ is polynomial and its inverse is defined whenever $\eta$ is invertible, Proposition~\ref{membership} gives $Z(L_\eta)\subseteq Z(\eta)$. At scalar size $t=0$, $v=r=1$ makes $\eta=1$, a contradiction. One may equivalently use its Schur polynomial linearization, whose determinant vanishes exactly when $\eta$ does. The seam replacement $q$ still omits $v$ and removes $V$.

Every chart contains $R_M$ by Proposition~\ref{mobius-coordinate}, so the common support of the fan and these flips is exactly $C$. Theorem~\ref{support-intersections} identifies the intersection of their images with $R_M$. The allowed algebra and the finite intersection both lie in the fan chart, where the map to the skew field is injective. This gives the equality in the normalized surface algebra. Restore the boundary free factor by Lemmas~\ref{gauge} and~\ref{freeproduct-intersections}. The full intersection is contained in the finite one and contains the allowed algebra by the Laurent and angle theorems. All presentations and formulas are defined over $\Q$, so finite-intersection scalar descent as in Theorem~\ref{polygon-theorem} gives the result over $\Q$.
\end{proof}

\section{The one-marked M\"obius obstruction}\label{sec-one-mark}
For the once-marked strip, lifted marked vertices are $p_i=g^ip_0$. Directed nontrivial relative classes $\gamma_j$ are indexed by $j\in\Z\setminus\{0\}$, with reversal $\gamma_{-j}$. The simple nontrivial classes have winding $\pm1,\pm2$. To see this, push a simple based loop into the interior. A noncontractible two-sided embedded circle lifts to an essential circle in the annular orientation cover and has winding $\pm2$ in the strip. A one-sided circle has doubled winding $\pm2$ and hence winding $\pm1$. The winding-two classes are the boundary orientations. The winding-one arc cuts the strip into a triangle, has no ordinary flip, and gives the unique ordinary triangulation. This is also the anti-self-folded case in \cite[Proposition 3.4 and Remark 6.3]{DP}.

Write $u=x_{\gamma_1}$, $v=x_{\gamma_{-1}}$, $z=x_{\gamma_2}$, $w=x_{\gamma_{-2}}$. The triangle relation and total angle are
\begin{equation}\label{one-mark-relations}
 z=uv^{-1}wv^{-1}u,\qquad
 T=v^{-1}uz^{-1}+u^{-1}zu^{-1}+z^{-1}uv^{-1}.
\end{equation}
The chart is the free group algebra on $u,v,w$.

\begin{theorem}\label{one-mark-theorem}
For this surface, the allowed curve, boundary-inverse, and total-angle algebra is strictly smaller than its unique ordinary chart. In particular $u^{-1}$ and $v^{-1}$ belong to the chart and do not belong to $\U_{M_1}$.
\end{theorem}
\begin{proof}
Over $\C(\varepsilon)$ put
\[
 Q=\begin{pmatrix}1&0\\0&-1\end{pmatrix},\quad
 S=\begin{pmatrix}0&1\\1&0\end{pmatrix},\quad
 C=\begin{pmatrix}0&\varepsilon+1\\\varepsilon&0\end{pmatrix},\quad
 K=\begin{pmatrix}S&0\\C&-S\end{pmatrix}.
\]
For $J=\left(\begin{smallmatrix}0&I_2\\-I_2&0\end{smallmatrix}\right)$, direct multiplication gives
$K^TJK=-J$ and $K^2=\left(\begin{smallmatrix}I&0\\Q&I\end{smallmatrix}\right)$.
Set $E_i=K^i\binom{I_2}{0}$. They are isotropic frames. Order the lifted vertices by reading even indices increasingly and then odd indices decreasingly; one deck step reverses this order. Signed pairings, with ordinary transpose in place of $*$, satisfy the polygon identities by the same calculation as Lemma~\ref{signed-pairings}, now with matrix blocks. The calculation is valid over the field $\C(\varepsilon)$; transverse frames and nondegenerate pairings provide the required inverses.

The power formulas give, for every integer $s$,
\begin{equation}\label{all-one-mark-curves}
 A_{2s}=|s|Q\ (s\neq0),\qquad
 A_{2s+1}=\begin{pmatrix}0&\varepsilon+s+1\\\varepsilon-s&0\end{pmatrix}.
\end{equation}
Every displayed matrix is invertible over $\C(\varepsilon)$ and $A_{-j}=A_j^T$. The fundamental values are
$u\mapsto C$, $v\mapsto C^T$, $z,w\mapsto Q$,
and satisfy \eqref{one-mark-relations}. Every finite consecutive union of lifted triangles is a polygon. Comparing its triangulation values and Laurent formula with these signed pairings proves that each curve evaluates to $A_j$. Thus all curves are invertible in the target, and the universal curve localization extends this chart homomorphism to the full surface algebra.

Writing $a=\varepsilon+1$, $b=\varepsilon$, one obtains
\[
 T\longmapsto Q\left(\frac ba+\frac ab-\frac1{ab}\right)=2Q,
\]
because $a-b=1$. All curve values in \eqref{all-one-mark-curves}, both boundary inverses, and this angle have polynomial entries in $\varepsilon$. In contrast
\[
 u^{-1}\longmapsto C^{-1}=\begin{pmatrix}0&\varepsilon^{-1}\\(\varepsilon+1)^{-1}&0\end{pmatrix}
\]
has a pole at zero. Therefore it is not in the allowed algebra. Reversal preserves that algebra: curves and boundary inverses are reversed, and each triangle angle is fixed by the triangle relation, so the total angle is fixed. Hence $v^{-1}$ is missing as well.
\end{proof}

\begin{proposition}[Residual extraction]\label{residual-extraction}
Use the polynomial core $P=\C\alg{u,v,w}$ and residuals at $(1,1,1)$. Then
\begin{equation}\label{extraction}
 \ell_w\ell_u r_u(uv^{-1}wv^{-1}u)=v^{-1}.
\end{equation}
The residual hull of the allowed algebra is the entire free group chart, and its additive saturation defect is infinite-dimensional.
\end{proposition}
\begin{proof}
With $z$ as in \eqref{one-mark-relations}, the product rules give
\[
 r_u z=uv^{-1}wv^{-1}+1,\qquad
 \ell_u r_u z=v^{-1}wv^{-1},\qquad
 \ell_w\ell_u r_u z=v^{-1}.
\]
These are coefficient residuals in the shifted variables, with operator composition from right to left. Also $z^{-1}=u^{-1}vw^{-1}vu^{-1}$ is an allowed boundary inverse and
$\ell_v\ell_v(z^{-1})=u^{-1}$,
by applying the left product rule twice. Since $w^{-1}$ is a boundary inverse, the hull contains all three inverse generators of the chart. Conversely that chart is residual-stable, being the inverse algebra of the three scalar pencils, and contains the allowed algebra. This proves the hull assertion, even though the allowed curve list is infinite: closure is defined by the same smallest-subalgebra construction. Theorem~\ref{infinite-defect} applies to the chart generated over the allowed algebra by these three inverses, and Theorem~\ref{one-mark-theorem} makes the defect nonzero.
\end{proof}
For the once-marked strip, the allowed algebra is strictly smaller than its residual hull, which equals the unique ordinary chart.

\section{Scope of the localization results}
The independent-holonomy theorems use a faithful universal localization of a free entry algebra. The free quotient calculation uses the new free coordinates in \eqref{free-quotient}. The involution calculation uses the pivot presentation \eqref{involution-chart}.

The polygon and M\"obius proofs use explicit coordinates for the normalized charts. To extend this approach to Berenstein--Retakh's general surface conjecture \cite[Conjecture 3.43]{BR}, one must determine the surface relation ideals and compute the tensor-layer corrections in Section~\ref{sec-descent}. The next step is to relate those correction classes to the allowed curve and angle generators. Decorated local-system coordinates give a geometric realization of the triangle and exchange identities \cite{KR,AGW}.

\appendix
\section{Finite exact support calculations}\label{appendix-effective}
\begin{proposition}\label{effective}
Suppose all coefficients are computable algebraic numbers and a common regular algebraic scalar point is supplied. Finite pencil presentations and their simple supports can be computed by finite exact operations. Membership of a supplied scalar-regular rational function in a specified finite pencil inverse algebra is decidable. Under the fir hypothesis of Theorem~\ref{relative-theorem}, relative divisors of scalar-invertible square matrices are likewise computable.
\end{proposition}
\begin{proof}
Build realizations by block formulas for sums, products, and inverses, or by the Schur construction in Lemma~\ref{linearization}. Reachable and observable spans stabilize in at most the state dimension; restriction to the reachable space and quotient by the unobservable space minimize the realization using exact linear algebra. Theorem~\ref{residual-theorem} supplies a finite residual generating set.

To compute composition factors of a coefficient tuple on $\C^n$, first span its generated unital algebra by successive products. This process stops in dimension at most $n^2$. If it is the full matrix algebra, the tuple is irreducible by Burnside's theorem. Otherwise a proper common invariant subspace exists. For each possible dimension, the finitely many Grassmannian charts give polynomial equations $A_iV=VB_i$ with a chart normalization for $V$. Solving these polynomial systems over the computable algebraic closure, by elimination, either finds such a subspace or proves its absence. If a complex solution exists, the algebraically closed coefficient field has a solution as well, by the Nullstellensatz. Recurse on the subspace and quotient to obtain a composition series and all multiplicities.

Discard zero tuples. Two irreducible tuples of equal size are simultaneously similar exactly when the linear system $A_iT=TB_i$ has a nonzero solution. Every nonzero intertwiner is invertible by irreducibility. Different sizes are distinct. This computes support containment, so Proposition~\ref{membership} decides membership in the specified inverse algebra. To compute the relative divisor, linearize the matrix and remove the base constituents with their multiplicities.
\end{proof}
For a finite list of holonomy entries, enumerate rational scalar tuples until all specified matrices and entries are invertible. The nonempty-open argument following Theorem~\ref{word-theorem} proves termination.

\section{A direct scalar regularity construction}\label{appendix-pingpong}
We give the scalar representation used in Theorem~\ref{word-theorem}. Choose $2r$ pairwise disjoint closed intervals $I_i^+,I_i^-$ in the real projective line. For each $i$ choose a hyperbolic projective transformation with attracting fixed point in $I_i^+$ and repelling fixed point in $I_i^-$. Sufficiently large positive powers give
\[
 T_i(\mathbb P^1\setminus I_i^-)\subset\operatorname{int} I_i^+,
 \qquad T_i^{-1}(\mathbb P^1\setminus I_i^+)\subset\operatorname{int} I_i^-.
\]
These inclusions follow by conjugating to a dilation and taking a large power. For a reduced nonempty word, choose a point outside the repelling interval of its final letter and outside the attracting interval of its first letter. Successive inclusions put its image inside the latter interval, so the word is nonidentity as a projective transformation. Lifts to $\GL_2(\R)$ thus have no nonidentity scalar words. A nonscalar operator has a vector $v$ with $Tv$ not proportional to $v$; a basis containing these vectors makes a specified off-diagonal coefficient nonzero after permutation. This proves the scalar nonvanishing used in the main argument for every finite rank.

\end{document}